\documentclass{article}
\usepackage{graphicx}
\usepackage{amsmath,amssymb,amsthm}
\usepackage{epsfig}
\usepackage{cite}
\usepackage{subfigure}
\usepackage{enumerate}
\usepackage{color}
\usepackage[a4paper, total={6.2in, 8.6in}]{geometry}
\allowdisplaybreaks
\theoremstyle{plain}
\newtheorem{thm}{Theorem}[section]
\newtheorem{lemm}[thm]{Lemma}
\theoremstyle{definition}

\newtheorem{defi}[thm]{Definition}
\newtheorem{cor}[thm]{Corollary}
\newtheorem{claim}[thm]{Claim}
\newtheorem{rmk}[thm]{Remark}
\newcommand{\defeq}{\stackrel{\text{def}}{=}}

\newcommand{\norm}[1]{\left\lVert#1\right\rVert}

\DeclareMathSymbol{\C}{\mathalpha}{AMSb}{"43}
\DeclareMathOperator{\supp}{supp}
\title{\bf Localization Results for Random Jacobi Operators}
\author{\textsc{ Valmir Bucaj}\footnote{The author was supported in part by NSF grant DMS--1361625.}}
\date{}
\usepackage{mathtools}

\DeclarePairedDelimiter\floor{\lfloor}{\rfloor}
\begin{document}

\nocite{*}

\title{The Kunz-Souillard Approach to Localization for Jacobi Operators }

\maketitle

\begin{abstract} \noindent In this paper we study spectral properties of Jacobi operators. In particular, we prove two main results: (1) that perturbing the diagonal coefficients of Jacobi operator, in an appropriate sense, results in exponential localization, and purely pure point spectrum with exponentially decaying eigenfunctions; and (2) we present examples of decaying {\it potentials} $b_n$ such that the corresponding Jacobi operators have purely pure point spectrum.
\end{abstract}
\vskip 12pt
\section{Introduction and Setting}
We will use the Kunz-Souillard approach to localization for random Schr\"odinger operators to prove that any Jacobi operator can be approximated by some {\it random} Jacobi operator, in operator norm, with purely pure point spectrum, and to also provide examples of Jacobi operators with decaying potentials having purely pure point spectrum.

\noindent The advantage of the Kunz-Souillard method is that it tackles localization directly, and you can add a background potential at no extra price. The shortcomings of the method are mainly because it applies only in one-dimension, and that it is known to work only for single-site distributions that are purely absolutely continuous, nevertheless, the conclusions are very strong. Whether this method can be extended to single-site distributions with a non-trivial singular part, still remains open.

\noindent Originally, the Kunz-Souillard work for Schr\"odinger operators was done in the discrete setting (see \cite{ks}). The analogue in the continuum setting was fully worked out by Damanik and Stolz (see \cite{ds}). 

\noindent Jacobi operators are important objects in mathematics. For one, they are a generalization of Schr\"odinger operators, which are central objects in quantum mechanics, also, the half line Jacobi operators with bounded coefficients correspond to compactly supported measures on the real line. Such correspondence can be established via orthogonal polynomials or the Borel transform of the measure. For a more elaborate discussion see \cite{oprl}.

\noindent On the other hand, the study of random Jacobi operators is of particular importance, since such operators model disordered media (e.g. amorphous solids). In some instances, as it is the case for crystals, the structure of the solid is completely regular; that is, the atoms are distributed periodically on some lattice. Then, mathematically, in such regular crystals, the total potential that a single particle (e.g. electron) at some position in $\mathbb{R}^d$, feels is periodic with respect to the lattice at hand. Schr\"odinger operators with periodic potentials are well understood, see for example \cite{rs} and \cite{eastham}.

\noindent However, as it is often the case in nature, if the positions of the atoms in the solid deviate from, say, a lattice in some highly non-regular way, or if the solid is some kind of mixture of various materials, then it is natural to view the potential that, say, a single particle feels at some position, as some random quantity. Mathematically, this can be studied via Jacobi operators with random potentials. So, understanding spectral properties of such operators is of great importance.\\

{\bf Set-up:} Suppose $r:\mathbb{R}\to \mathbb{R}_{\geq 0}$  is bounded, measurable, and compactly supported with $\norm{r}_1=1.$ Let $c\in\ell^{\infty}(\mathbb{Z})$.  Define a measure $\mu_n$ on $\mathbb{R}$ via
 $d\,\mu_n(E)=r_n(E)dE,$ where $r_n(x)=d_n^{-1}r\left(d_n^{-1}x\right),$ and $d_n$ is some fixed sequence.  
 Let

\begin{eqnarray*}
M&=&\sup\{|E|:E\in \supp(r)\}\\
M_n&=&\sup\{|E|: E\in \supp(r_n)\}\\
I_n&=&[c(n)-M_n, c(n)+M_n]\\
\Omega&=&\prod_{n\in\mathbb{Z}}I_n\\
d\mu(x)&=&\prod_{n\in \mathbb{Z}}r_n(x_n-c(n))dx_n.
\end{eqnarray*}

We wish to point out that $r$ quantifies the deviation of our random potential from the background potential $c$. In the second situation we will consider, the sequence $d_n$ will serve as a damping parameter that we will use to force decay of the random potential.

\noindent Next, we define $b_{\omega}(n)=\omega(n)$ for each $\omega\in \Omega.$ Notice, that each $b_{\omega}(n)$ is the sum of a random $i.i.d$ with distribution $\mu_n$ and some fixed background {\it potential} $c(n).$

\noindent With this notation, we define a one parameter family of Jacobi operators, $J_{\omega}$, on $\ell^2(\mathbb{Z})$ as follows

\begin{equation}\label{eq00}\left(J_{\omega}\phi\right)(n)=a(n)\phi(n+1)+a(n-1)\phi(n-1)+b_{\omega}(n)\phi(n),\end{equation}
where $a\in\ell^{\infty}(\mathbb{Z})$ with $a(n)\geq \delta>0$ for all $n\in \mathbb{Z}.$ \\

In general, if one assumes that supp $ r$ contains more than one element-by construction, this is the case for us-the resulting family $\{J_\omega\}_{\omega\in\Omega}$ of operators, with $a(n)\equiv 1$ and $d_n\equiv 1$,  is referred to as the {\it Anderson model}. The, simplest non-trivial case, where supp $r$ contains precisely two elements is known as the {\it Bernoulli-Anderson model}. It is well known, that the spectrum of the Anderson model has a simple description, namely we have $$\Sigma_{J_\omega}=\Sigma\defeq\big[-2\norm{a}_\infty,2\norm{a}_\infty\big]+ \supp( r)\defeq\big\{a+b:a\in\big[-2\norm{a}_\infty,2\norm{a}_\infty\big], b\in\supp(r)\big\},$$ for $\mu-$ almost every $\omega\in\Omega.$ This description shows that the spectrum of an Anderson model will always be given by a finite union of compact intervals. For a more extensive discussion see \cite{gs}.

\noindent One interesting property to study for the Anderson model is the phenomenon of {\it localization}. There are typically two separate statements referring to localization: a spectral statement and a dynamical one. {\it Spectral Anderson localization} asserts that the operators $J_\omega$ almost surely have pure point spectrum, with exponentially decaying eigenfunctions. On the other hand, different notions of {\it dynamical Anderson localization} have been used in literature. However, in essence, dynamical localization refers to an absence of transport in a random medium. This is typically quantified via (almost-sure) bounds on the moments of wave packets such as
$$\sup_{t}\sum_{n\in\mathbb{Z}}|n|^p\left|\langle{\delta_n,e^{-itJ_\omega}\delta_0\rangle}\right|^2<\infty,$$ for all $p>0.$ In some instances, one can prove stronger statements, such as replacing the {\it almost sure} condition by an expectation $\mathbb{E}(\cdot),$ as is the case via the Kunz-Souillard approach to localization in dimension one. For a more elaborate discussion of this method in the case of Schr\"odinger operators, see \cite{cfks}. 

In an appropriate formulation, it is known that dynamical localization implies spectral localization, while the converse is not true. For example, the so called random dimer model serves as a counterexample to this implication (see \cite{JS} and \cite {JSS} for a more elaborate description). One typically needs `` spectral localization +$\epsilon$" to imply dynamical localization in some suitable formulation. This relationship was studied by del Rio, Jitomirskaya, Last, and Simon in \cite{rjls}.

There are different approaches to localization: {\it Spectral averaging} can be used to study spectral localization; One can also study both spectral and dynamical localization via methods such as, {\it multi-scale analysis}, {\it fractional moments method} and also, which is what we do in this paper, the {\it Kunz-Souillard} method. Each method has its advantages and disadvantages. On one hand, some are broader in generality, but the conclusions one obtains are typically weaker. On the other hand, as is the case with the Kunz-Souillard approach, the scope of generality is narrow, however, the conclusions one draws are very rich.

\section{Main Results}
Our main goal is to prove the following two theorems which establish spectral localization for the family $(J_\omega)_{\omega\in\Omega}$ under suitable conditions.

\begin{thm}\label{thm2.0} For all $a_n, b_n\in \mathbb{R}$ bounded, with $a_n\geq \delta>0$, and for every $\epsilon>0$, there exist $\tilde a_n, \tilde b_n$, with $\norm{\tilde a-a}_\infty<\epsilon$ and $\norm{\tilde b-b}_\infty<\epsilon$, such that the Jacobi operator, $\widetilde{J}\defeq \widetilde{J}(\tilde a,\tilde b)$, has purely pure point spectrum with exponentially decaying eigenfunctions.
\end{thm}

\begin{thm}\label{thmm2.2}
With the same notation as above, if $d_n$ is a fixed sequence with $0\leq d_n\leq1$ and $d_n\geq C|n|^{-\zeta}$ for $\zeta<\frac{1}{2}$, then for $\mu$-almost every $\omega$, the Jacobi operator $J_\omega$ has purely pure point spectrum.
\end{thm}

The key ingredients of the proofs of Theorems $\ref{thm2.0}$ and $\ref{thmm2.2}$ are given by the following theorems, which are important in their own right, since they establish dynamical localization. 

\begin{thm} \label{thm2.1}With $\Omega$, $\mu$, and $J_{\omega}$ as above, and $d_n=1$ for all $n$, there exist constants $C,\gamma\in(0,\infty)$ such that 
$$\int_{\Omega}\left(\sup_{t\in\mathbb{R}}\left|\langle{\delta_m, e^{-itJ_{\omega}}\delta_0}\rangle\right|\right)d\mu(\omega)\leq Ce^{-\gamma|m|},$$ for all $m\in \mathbb{Z}.$
\end{thm}

Actually, we can loosen the condition on the sequence $d_n$; that is, the statement holds true as long as $d_n\in \ell^{\infty}(\mathbb{R})$ is positive and uniformly bounded away from zero.

For more pleasant exposition let $$a(m,n)=\int_{\Omega}\left(\sup_{t\in\mathbb{R}}\left|\langle{\delta_m, e^{-itJ_{\omega}}\delta_n}\rangle\right|\right)d\mu(\omega).$$

\begin{rmk}\label{rmk01}
We wish to point out that in a similar way one shows that \begin{equation}a(m,n)\leq Ce^{-\gamma|m-n|}.\end{equation} For simplicity, we only work out the case $n=0.$
\end{rmk}

\begin{thm}\label{thm2}
If there are constants $C,\gamma\in (0,\infty)$ such that $$\max_{n\in\{0,1\}}a(m,n)\leq Ce^{-\gamma|m|},$$ then for $\mu-$almost every $\omega\in \Omega$, $J_\omega$(this is as in $(\ref{eq00})$) has pure point spectrum with exponentially decaying eigenfunctions. More precisely, these eigenfunctions obey estimates of the form $$|u(m)|\leq C_{\omega, \epsilon, u}e^{-(\gamma-\epsilon)|m|},$$ for small enough $\epsilon\in(0,\gamma).$
\end{thm}

\begin{proof}
 This is proved in almost identical way as in the case for random Schr\"{o}dinger operators, so we direct the reader to \cite{dd} or \cite{cfks}.

\end{proof}
Even if we do not insist on exponential bounds for $\displaystyle\max_{n\in\{0,1\}}a(m,n)$, we still obtain pure point spectrum, but we no longer get exponentially decaying eigenfunctions. We make this statement precise in the following two theorems. 

\begin{thm}\label{thmm2.6} Let  $d_n$ be a fixed sequence with $0\leq d_n\leq1$ and $d_n\geq C|n|^{-\zeta}$ for $\zeta<\frac{1}{2}$ and some constant $C>0$. With $\Omega$, $\mu$, and $J_{\omega}$ as above, there exist constants $C'>0$ and $\gamma''>0,$ such that 
$$\int_{\Omega}\left(\sup_{t\in\mathbb{R}}\left|\langle{\delta_m, e^{-itJ_{\omega}}\delta_0}\rangle\right|\right)d\mu(\omega)\leq C' |m|^{\zeta/2}\exp\left(-\gamma''|m|^{1-2\zeta}\right).$$ 
\end{thm}
\begin{rmk}
As in Remark \ref{rmk01}, we only work out the proof for $a(m,0)$, since the other cases are completely analogous.
\end{rmk}
\begin{thm}\label{thmm2.5}
If there exist constants $C''>0$ and $\tau>\frac{3}{2}$, such that $$\max_{n\in\{0,1\}}a(m,n)\leq \frac{C''}{m^{\tau}},$$  then for $\mu$-almost every $\omega\in \Omega$, the Jacobi operator $J_\omega$, has purely pure point spectrum.
\end{thm}

\begin{proof}
Let us define $$a(m,n,\omega)=\sup_{t\in\mathbb{R}}\left|\langle{\delta_m,e^{-itJ_\omega}\delta_n\rangle}\right|,$$ so that we have $$a(m,n)=\int_{\Omega}a(m,n,\omega)d\mu(\omega).$$
Let $\frac{1}{2}<\beta<\tau-1$ be given, and consider the set 

$$S_{\beta,m,n}=\Big\{\omega\in\Omega: a(m,n,\omega)>\frac{1}{m^\beta}\Big\}.$$

\noindent Then $$a(m,n)\geq \frac{1}{m^\beta}\mu\left(S_{\beta,m,n}\right),$$ for all $m,n\in\mathbb{Z}.$ So, by the above observation and the hypothesis, for all $m$, and $n=0,1$ we get \begin{equation}\label{eqeq1}\mu\left(S_{\beta,m,n}\right)\leq m^\beta a(m,n)\leq\frac{C''}{m^{\tau-\beta}}.\end{equation}

 \noindent Thus, since by our choice of $\beta$ we have $\tau-\beta>1$, by comparison test, from $(\ref{eqeq1})$, we get $$\sum_{m\in\mathbb{Z}}\mu\left(S_{\beta,m,n}\right)<\infty,$$ for $n=0,1.$ As a result, by Borel-Cantelli lemma, we have

$$\mu\left( \Big\{\omega\in\Omega: a(m,n,\omega)>\frac{1}{m^\beta},\, \mbox{for infinitely many $m$}\Big\}  \right)=0.$$

\noindent Let $$\Omega_0=\Big\{\omega\in\Omega: a(m,n,\omega)\leq\frac{1}{m^\beta}, \,\mbox{for all but finitely many $m$}\Big\},$$ for $n=0,1$, with $\mu\left(\Omega_0\right)=1.$ Then, it follows that for all $\omega\in\Omega_0$ we have

$$C_{\omega,\beta}\defeq\sup_{n=0,1,m\in\mathbb{Z}} a(m,n,\omega)m^\beta<\infty.$$

As a consequence, we get $$a(m,n,\omega)\leq \frac{C_{\omega,\beta}}{m^\beta}.$$ In particular, for each fixed $M>0$, we have
\begin{equation}\label{eqeq2}\sum_{|m|\geq M}\left|\langle{\delta_m, e^{-itJ_\omega}\delta_n\rangle}\right|^2\leq \sum_{|m|\geq M}\left(a(m,n,\omega)\right)^2\leq\sum_{|m|\geq M}C_{\omega,\beta}^2\frac{1}{m^{2\beta}},\end{equation} for $n=0,1.$

\noindent Since, by assumption, we have $\beta>1/2$, it follows that the series in $(\ref{eqeq2})$ goes to zero as $M\to\infty$.
 In particular, for every $\epsilon>0$, there is some $M>0$ such that 
$$\sum_{|m|\geq M}\left|\langle{\delta_m, e^{-itJ_\omega}\delta_n\rangle}\right|^2<\epsilon,$$ for every $t\in\mathbb{R},$ and $n=0,1.$

\noindent Thus, by RAGE theorem, it follows that the spectral measures $\mu_{\delta_0}^{J_\omega}$ and $\mu_{\delta_1}^{J_\omega}$ are pure point measures. On the other hand, since the pair $\{\delta_0,\delta_1\}$ is a spectral basis for the operator $J_\omega$, it follows that all spectral measures of $J_\omega$ are pure point measures. So, in conclusion, for each $\omega\in \Omega_0$, the Jacobi operator $J_\omega$, has purely pure point spectrum.

\end{proof}


\subsection{Proof of Theorem \ref{thm2.0}}

\begin{proof}
Let $J\defeq J(a_n,b_n)$ be a given Jacobi operator, where $a_n,b_n$ are as in the statement of the theorem. Given $\epsilon >0$ we will construct $\widetilde{J}=\widetilde{J}(\tilde a,\tilde b)$ as follows. We pick $\tilde a\defeq a$, and $\tilde b(n)\defeq b_\omega(n)$, where $b_\omega(n)$ is as above, with $c(n)$ replaced by $b(n)$ and $M<\epsilon$. 
Then, clearly$ \norm{\tilde a-a}_\infty<\epsilon$ and $\norm{\tilde b-b}_\infty<\epsilon$. Then, by Theorems \ref{thm2.1} and \ref{thm2}, it follows that $\widetilde{J}$ has purely pure point spectrum. 

\end{proof}
\begin{rmk} We have shown the much stronger statement; that is, we showed that there exist and uncountable family of operators with the desired property.
\end{rmk}

\subsection{Proof of Theorem \ref{thmm2.2}}

\begin{proof}
This is an immediate consequence of Theorems \ref{thmm2.5}. More specifically, we claim that for large enough $m$ and some $\tau>3/2,$ we have

\begin{equation}\label{eqqq}|m|^{\zeta/2}\exp\left(-\gamma'' \left|m \right|^{1-2\zeta}\right)\leq \frac{1}{m^{\tau}}.\end{equation}

\noindent A quick calculation shows that $$\lim_{m\to\infty}|m|^{\zeta/2+\tau}\exp\left(-\gamma'' \left|m\right|^{1-2\zeta}\right)= 0,$$ which, in turn, implies $(\ref{eqqq})$. Then, this observation and Theorem \ref{thmm2.6} imply that for $n=0,1,$ we have $$a(m,n)\leq\frac{C''}{m^{\tau}}.$$ Thus, the result follows from Theorem \ref{thmm2.5}. 

\end{proof}







\section{Preparatory Work}
We turn to the task of proving Theorems $\ref{thm2.1}$ and $\ref{thmm2.6}$.

Given $L\in \mathbb{Z}_+$, denote by $J_{\omega}^{(L)}$ the restriction of $J_{\omega}$ to $\ell^2(-L,\dots, L)$, and let 
$$a_L(m,n)=\int_{\Omega}\left(\sup_{t\in\mathbb{R}}\left|\langle{\delta_m, e^{-itJ_{\omega}^{(L)}}\delta_n}\rangle\right|\right)d\mu(\omega).$$ That is, 
\begin{equation*}
  J_{\omega}^{(L)}=\left(
  \begin{array}{cccccc} 
    b_{\omega}(-L)& a(-L)& &&&  0 \\ 
    a(-L)& b_{\omega}(-L+1)&&& &  0\\ 
 0&a(-L+1)&&&&\\
    \vdots&  &&\ddots&a(L-2)  & \vdots \\ 
    &&&&b_{\omega}(L-1)&a(L-1)\\
    0& &&\ldots&a(L-1) & b_{\omega}(L) \\ 
  \end{array}
  \right).
\end{equation*}

Let $\{E_{\omega}^{L,k}\}_k,$ and $\{\varphi_{\omega}^{L,k}\}_k$ be the eigenvalues and the corresponding normalized eigenfunctions of $J_{\omega}^{(L)}$, respectively. Define,
$$\rho_L(m,n)=\int_{\Omega}\left(\sum_{k}\left|\langle{\delta_m,\varphi_{\omega}^{L,k}}\rangle\right|\left|\langle{ \delta_n,\varphi_{\omega}^{L,k}}\rangle\right|\right)d\mu(\omega),$$
and notice that this is a $(2L+1)$ fold integral, since  $J_{\omega}^{(L)}$ depends only on the entries $\omega_{-L},\dots\omega_{L}.$\\

The following two lemmas are easy to prove, for a discussion see \cite[pp. 192-193]{cfks}. However, for completeness and reader's convenience, we include the brief arguments here.

\begin{lemm} For $m,n\in \mathbb{Z}$ we have $$a(m,n)\leq\liminf_{L\to\infty}a_L(m,n).$$
\end{lemm}

\begin{proof}
First, regarding $J_{\omega}^{(L)}$ as an operator in $\ell^2(\mathbb{Z})$ , in the natural way, we observe that $J_{\omega}^{(L)}$ converges strongly to $J_{\omega}.$ As a consequence, $e^{-itJ_{\omega}^{(L)}}$ converges strongly to $e^{-itJ_{\omega}}$, for each $t\in \mathbb{R},$ and every $\omega.$ As a result, we also have
$$\lim_{L\to \infty}\left|\langle{\delta_m, e^{-itJ_{\omega}^{(L)}}\delta_n}\rangle\right|=\left|\langle{\delta_m, e^{-itJ_{\omega}}\delta_n}\rangle\right|.$$ Next, for each $t\in \mathbb{R}$, we have 
$$\left|\langle{\delta_m, e^{-itJ_{\omega}^{(L)}}\delta_n}\rangle\right|\leq\sup_{t'\in\mathbb{R}}\left|\langle{\delta_m, e^{-it'J_{\omega}^{(L)}}\delta_n}\rangle\right|.$$ Taking lim inf of both sides we obtain:

$$\left|\langle{\delta_m, e^{-itJ_{\omega}}\delta_n}\rangle\right|=\lim_{L\to \infty}\left|\langle{\delta_m, e^{-itJ_{\omega}^{(L)}}\delta_n}\rangle\right|\leq\liminf_{L\to\infty}\sup_{t'\in\mathbb{R}}\left|\langle{\delta_m, e^{-it'J_{\omega}^{(L)}}\delta_n}\rangle\right|.$$ Hence, 
$$\sup_{t\in\mathbb{R}}\left|\langle{\delta_m, e^{-itJ_{\omega}}\delta_n}\rangle\right|\leq\liminf_{L\to\infty}\sup_{t\in\mathbb{R}}\left|\langle{\delta_m, e^{-itJ_{\omega}^{(L)}}\delta_n}\rangle\right|.$$
The result follows by an application of Fatou's lemma.
\end{proof}

\begin{lemm}For $L\in\mathbb{Z}_+$, and $m,n\in\mathbb{Z}$ we have $$a_L(m,n)\leq\rho_L(m,n).$$
\end{lemm}

\begin{proof}
We have

\begin{align*}
a_L(m,n)&=\int_{\Omega}\left(\sup_{t\in\mathbb{R}}\left|\langle{\delta_m, e^{-itJ_{\omega}^{(L)}}\delta_n}\rangle\right|\right)d\mu(\omega)\\
&=\int_{\Omega}\left(\sup_{t\in\mathbb{R}}\left|\langle{\delta_m, e^{-itJ_{\omega}^{(L)}}\sum_k\langle{\delta_n,\varphi_{\omega}^{L,k}}\rangle}\varphi_{\omega}^{L,k}\rangle\right|\right)d\mu(\omega)\\
&\leq\int_{\Omega}\left(\sup_{t\in\mathbb{R}}\sum_k \left|\langle{\delta_m, e^{-itE_{\omega}^{L,k}}\langle{\delta_n,\varphi_{\omega}^{L,k}}\rangle}\varphi_{\omega}^{L,k}\rangle\right|\right)d\mu(\omega)\\
&=\int_{\Omega}\left(\sup_{t\in\mathbb{R}}\sum_k \left|\langle{\delta_m,\varphi_{\omega}^{L,k}}\rangle\right|\left|\langle{\delta_n,\varphi_{\omega}^{L,k}}\rangle\right|\right)d\mu(\omega)\\
&=\rho_L(m,n).
\end{align*}
\end{proof}

Put $$\Sigma_0=\left[-2\norm{a}_{\infty}-M-\norm{c}_{\infty}, 2\norm{a}_{\infty}+M+\norm{c}_{\infty}\right].$$ Notice that $\Sigma_0$ contains the spectrum of both $J_{\omega},$ and $J_{\omega}^{(L)}$. Now, in the spirit of \cite{ks}, we define a family of operators appropriate for our setting.

\begin{defi}
For $E\in \mathbb{R}$, define the operators $U, S_E^{(n)}, T_E^{(n)}$ on $L^p(\mathbb{R})$ by:
$$\left(Uf\right)(x)=|x|^{-1}f(x^{-1}).$$
\begin{displaymath}
   \left(S_E^{(n)}f\right)(x) = \left\{
     \begin{array}{lr}
       \displaystyle a_n\int r_n(E-a_nx-a_{n-1}y^{-1})f(y)dy& , n<0\\
      \displaystyle a_0\int r_0(E-a_0x-a_{-1}y^{-1})f(y)dy& , n=0\\
\displaystyle a_{n-1}\int r_n(E-a_{n-1}x-a_ny^{-1})f(y)dy&,n>0
     \end{array}
   \right.
\end{displaymath} 

and 
\begin{eqnarray*}
\left(T_E^{(n)}f\right)(x)&=&\sqrt{a_{n-1}a_n}\int r_n(E-a_{n-1}x-a_ny^{-1})|y|^{-1}f(y)dy, \,n>0.\\
r_{k;E}^{(n)}(x)&=&r_k(E-a_{n-1}x)\\
\end{eqnarray*}
We also need to define the following: 

\begin{eqnarray*}
S_{E;m}^{(n)}&=&S_{E-c(m)}^{(n)};\\
T_{E;m}^{(n)}&=&T_{E-c(m)}^{(n)};\\
r_{k;E;m}^{(n)}&=&r_{k;E-c(m)}^{(n)}.
\end{eqnarray*}

We wish to point out that $U$ is a unitary operator on $L^2(\mathbb{R}).$
\end{defi}

From now on, we will drop the subscript $\omega$ on the sequence $b$ (i.e. $b_n=b_\omega(n)=\omega(n)$), this should cause no confusion and should be clear from the context. 

We want to compute the following:

\begin{align}\label{eq0}\nonumber \rho_L(m,0)&=\int_{\Omega}\left(\sum_{k}\left|\langle{\delta_m,\varphi_{\bar{b}}^{L,k}}\rangle\right|\left|\langle{ \delta_0,\varphi_{\bar b}^{L,k}}\rangle\right|\right)d\mu(\omega)\\
&=\int\dots\int\left(\sum_{k}\left|\langle{\delta_m,\varphi_{\bar{b}}^{L,k}}\rangle\right|\left|\langle{ \delta_0,\varphi_{\bar b}^{L,k}}\rangle\right|\right)\prod_{n=-L}^Lr_n(b_n-c_n)db_{-L}\dots db_L,
\end{align}

where $\bar b=(b_{-L},\dots, b_L).$
 Let $\{E_{\bar b}^{L,k}\}_{-L\leq k\leq L}$ and $\{\varphi_{\bar b}^{L,k}\}$ be the eigenvalues and the corresponding normalized eigenvectors of

\begin{equation*}
  J_{\omega}^{(L)}=\left(
  \begin{array}{cccccc} 
    b_{-L}& a_{-L}& &&&  0 \\ 
    a_{-L}& b_{-L+1}&&& &  \\ 
 0&a_{-L+1}&&&&\\
    \vdots&  &&\ddots&a_{L-2}  & \vdots \\ 
    &&&&b_{L-1}&a_{L-1}\\
    0& &&\ldots&a_{L-1} & b_L \\ 
  \end{array}
  \right).
\end{equation*}

Let $E$ be $E_{\bar b}^{L,k}$ and $u$ be $\varphi_{\bar b}^{L,k}$, then we have 
\begin{equation}\label{eq1}a_nu_{n+1}+a_{n-1}u_{n-1}+b_nu_n=Eu_n,\end{equation} for $-L\leq n\leq L$, where $u_{-L-1}=u_{L+1}=0.$

Rewriting $(\ref{eq1})$ we get: 
\begin{equation}\label{eq2} b_n=E-a_n\frac{u_{n+1}}{u_n}-a_{n-1}\frac{u_{n-1}}{u_n}\end{equation}
Let

\begin{displaymath}
   x_n= \left\{
     \begin{array}{lr}
       \displaystyle \frac{\varphi_{\bar b}^{L,k}(n+1)}{\varphi_{\bar b}^{L,k}(n)}& , n<0\\
\displaystyle \frac{\varphi_{\bar b}^{L,k}(n-1)}{\varphi_{\bar b}^{L,k}(n)} &,n>0
     \end{array}
   \right.
\end{displaymath} 
so that 

\begin{displaymath}
   b_n= \left\{
     \begin{array}{lr}
       \displaystyle E-a_{n-1}x_{n-1}^{-1}-a_nx_n ,&  n<0\\
\displaystyle E-a_{-1}x_{-1}^{-1}-a_0x_1^{-1}, &n=0\\
\displaystyle E-a_nx_{n+1}^{-1}-a_{n-1}x_n,&n>0
     \end{array}
   \right.
\end{displaymath}  with the convention $x_{-L-1}^{-1}=x_{L+1}^{-1}=0.$ 

This motivates the following change of variables $$F_L:(x_{-L},\dots, x_{-1}, E, x_1,\dots, x_L)\mapsto(b_{-L}, \dots,b_0,\dots, b_L).$$

The next step is to rewrite $(\ref{eq0})$ using this change of variables. In order to do so, we need to compute the determinant of the Jacobian of this change of variables. 

Observe that: 
$\frac{ \partial b_n}{\partial E}=1$, for all n; $\frac{\partial b_n}{\partial x_n}=-a_n$, for $n<0$; $\frac{\partial b_n}{\partial x_n}=-a_{n-1}$, for $n>0;$ $\frac{\partial b_n}{\partial x_{n-1}}=a_{n-1}x_{n-1}^{-2},$ for $n\leq 0$; $\frac{\partial b_n}{\partial x_{n+1}}=a_nx_{n+1}^{-2}$, for $n\geq 0;$ and $\frac{\partial b_n}{\partial x_m}=0$, for all other $m,n$. 

Thus, the corresponding matrix of $F_L$ is:

\begin{equation*}
  \left(
  \begin{array}{cccccccccccc} 
    -a_{-L}& a_{-L}x_{-L}^{-2}& &&&  &&&&& \\ 
    & -a_{-L+1}&a_{-L+1}x_{-L+1}^{-2}&&&&& &&&  \\ 
 &&-a_{-L+2}&&&&&&&&\\
    &  &&\ddots&\ddots  & &&&&& \\ 
&&&&-a_{-1}&a_{-1}x_{-1}^{-2}&&&&&\\
1&1&\ldots&&1&1&1&\ldots&&1&1\\
&&&&&a_0x_1^{-2}&-a_0&&&&\\
    &&&&&&a_1x_2^{-2}&-a_{1}&&&\\
    & &&&& &&\ddots&\ddots&&\\ 
&&&&&&&&a_{L-2}x_{L-1}^{-2}&-a_{L-2}&\\
&&&&&&&&&a_{L-1}x_L^{-2}&-a_{L-1}\\
  \end{array}
  \right).
\end{equation*}

We claim that 
\begin{align}\label{eq3}\nonumber\det F_L&=\left(\prod_{n=-L}^{L-1}a_n\right)\Big(1+x_1^{-2}\{1+x_2^{-2}\{1+\dots x_{L-1}^{-2}\{1+x_L^{-2}\}\dots\}\}\Big.\\
&\qquad+\Big.x_{-1}^{-2}\{1+x_{-2}^{-2}\{1+\dots x_{-L+1}^{-2}\{1+x_{-L}^{-2}\}\dots\}\}\Big)\\ \nonumber
&=\left(\prod_{n=-L}^{L-1}a_n\right)\left(\varphi_{\bar b}^{L,k}(0)\right)^{-2}.\\ \nonumber
\end{align}

We prove this by induction on $L$. For $L=1$ it is clear. Now, suppose that $(\ref{eq3})$ holds for some $L$. Consider the determinant of matrix of $F_{L+1}$:

\begin{equation*}
 \left(
  \begin{array}{ccccccccccc} 
    -a_{-L-1}& a_{-L-1}x_{-L-1}^{-2}& &&&  &&&&& \\ 
    & -a_{-L}&a_{-L}x_{-L}^{-2}&&&&& &&&  \\ 
 &&-a_{-L+1}&&&&&&&&\\
    &  &&\ddots&\ddots  & &&&&& \\ 
&&&&-a_{-1}&a_{-1}x_{-1}^{-2}&&&&&\\
1&1&\ldots&&1&1&1&\ldots&&1&1\\
&&&&&a_0x_1^{-2}&-a_0&&&&\\
    &&&&&&a_1x_2^{-2}&-a_{1}&&&\\
    & &&&& &&\ddots&\ddots&&\\ 
&&&&&&&&a_{L-1}x_{L}^{-2}&-a_{L-1}&\\
&&&&&&&&&a_{L}x_{L+1}^{-2}&-a_{L}\\
  \end{array}
  \right).
\end{equation*}

Expanding along the first column we get: 

\begin{equation*}
(-a_{-L-1}) \det\left(
  \begin{array}{ccccccccccc} 
     -a_{-L}&a_{-L}x_{-L}^{-2}&&&&& &&&  \\ 
 &-a_{-L+1}&&&&&&&&\\
      &&\ddots&\ddots  & &&&&& \\ 
&&&-a_{-1}&a_{-1}x_{-1}^{-2}&&&&&\\
1&\ldots&&1&1&1&\ldots&&1&1\\
&&&&a_0x_1^{-2}&-a_0&&&&\\
    &&&&&a_1x_2^{-2}&-a_{1}&&&\\
     &&&& &&\ddots&\ddots&&\\ 
&&&&&&&a_{L-1}x_{L}^{-2}&-a_{L-1}&\\
&&&&&&&&a_{L}x_{L+1}^{-2}&-a_{L}\\
  \end{array}
  \right)
\end{equation*}

+
\begin{equation*}
(-1)^{L+1} \det\left(
  \begin{array}{ccccccccccc} 
     a_{-L-1}x_{-L-1}^{-2}& &&&  &&&&& \\ 
     -a_{-L}&a_{-L}x_{-L}^{-2}&&&&& &&&  \\ 
 &-a_{-L+1}&&&&&&&&\\
      &&\ddots&\ddots  & &&&&& \\ 
&&&-a_{-1}&a_{-1}x_{-1}^{-2}&&&&&\\
&&&&a_0x_1^{-2}&-a_0&&&&\\
    &&&&&a_1x_2^{-2}&-a_{1}&&&\\
     &&&& &&\ddots&\ddots&&\\ 
&&&&&&&a_{L-1}x_{L}^{-2}&-a_{L-1}&\\
&&&&&&&&a_{L}x_{L+1}^{-2}&-a_{L}\\
  \end{array}
  \right).
\end{equation*}

Note that the second matrix is lower-triangular, so expanding along the first row, repeatedly, we eventually will get: $$\left(\prod_{n=-L-1}^La_n\right) x_{-L-1}^{-2}x_{-L}^{-2}\dots x_{-1}^{-2}.$$
Expanding the first determinant along the last column we  get: 

\begin{equation*}
(-1)^L \det\left(
  \begin{array}{ccccccccccc} 
     -a_{-L}&a_{-L}x_{-L}^{-2}&&&&& &&  \\ 
 &-a_{-L+1}&&&&&&&\\
      &&\ddots&\ddots  & &&&& \\ 
&&&-a_{-1}&a_{-1}x_{-1}^{-2}&&&&\\
&&&&a_0x_1^{-2}&-a_0&&&\\
    &&&&&a_1x_2^{-2}&-a_{1}&&\\
     &&&& &&\ddots&\ddots&\\ 
&&&&&&&a_{L-1}x_{L}^{-2}&-a_{L-1}\\
&&&&&&&&a_{L}x_{L+1}^{-2}\\
  \end{array}
  \right)
\end{equation*}

+
\begin{equation*}
(-a_L) \det\left(
  \begin{array}{ccccccccccc} 
     -a_{-L}&a_{-L}x_{-L}^{-2}&&&&& && \\ 
 &-a_{-L+1}&&&&&&&\\
      &&\ddots&\ddots  & &&&& \\ 
&&&-a_{-1}&a_{-1}x_{-1}^{-2}&&&&\\
1&\ldots&&1&1&1&\ldots&&1\\
&&&&a_0x_1^{-2}&-a_0&&&\\
    &&&&&a_1x_2^{-2}&-a_{1}&&\\
     &&&& &&\ddots&\ddots&\\ 
&&&&&&&a_{L-1}x_{L}^{-2}&-a_{L-1}\\
  \end{array}
  \right).
\end{equation*}

As before, computing the fist determinant by expanding along the first columns, repeatedly, we eventually get:

$$\left(\prod_{n=-L}^{L}a_n\right)x_1^{-2}x_2^{-2}\dots x_L^{-2}x_{L+1}^{-2}.$$
Combining all of these, and noting that the last determinant is simply $\det F_L$ we get:

\begin{align*}\det F_{L+1}&=(-a_{-L-1})\left(\left(\prod_{n=-L}^{L}a_n\right)x_1^{-2}x_2^{-2}\dots x_L^{-2}x_{L+1}^{-2}+(-a_L)\det F_L\right)+\left(\prod_{n=-L-1}^La_n\right) x_{-L-1}^{-2}x_{-L}^{-2}\dots x_{-1}^{-2}\\
&=a_{-L-1}a_L\det F_L+\left(\prod_{n=-L-1}^{L}a_n\right)x_1^{-2}x_2^{-2}\dots x_L^{-2}x_{L+1}^{-2}
+\left(\prod_{n=-L-1}^La_n\right) x_{-L-1}^{-2}x_{-L}^{-2}\dots x_{-1}^{-2}\\
&=a_{-L-1}a_L\prod_{n=-L}^{L-1}a_n\Big(1+x_1^{-2}\{1+x_2^{-2}\{1+\dots x_{L-1}^{-2}\{1+x_L^{-2}\}\dots\}\}\Big.\\
&\qquad+\Big.x_{-1}^{-2}\{1+x_{-2}^{-2}\{1+\dots\{x_{-L+1}^{-2}\{1+x_{-L}^{-2}\}\dots\}\}\Big)\\
&+\prod_{n=-L-1}^{L}a_nx_1^{-2}x_2^{-2}\dots x_L^{-2}x_{L+1}^{-2}+\prod_{n=-L-1}^La_n x_{-L-1}^{-2}x_{-L}^{-2}\dots x_{-1}^{-2}\\
&=\prod_{n=-L-1}^La_n\Big(1+x_1^{-2}\{1+x_2^{-2}\{1+\dots x_{L-1}^{-2}\{1+x_L^{-2}\}\dots\}\}\Big.\\
&\qquad+\Big.x_{-1}^{-2}\{1+x_{-2}^{-2}\{1+\dots\{x_{-L+1}^{-2}\{1+x_{-L}^{-2}\}\dots\}\}+x_1^{-2}x_2^{-2}\dots x_L^{-2}x_{L+1}^{-2}+ x_{-L-1}^{-2}x_{-L}^{-2}\dots x_{-1}^{-2}\Big)\\
&=\prod_{n=-L-1}^La_n\Big(1+x_1^{-2}\{1+x_2^{-2}\{1+\dots x_{L}^{-2}\{1+x_{L+1}^{-2}\}\dots\}\}\Big.\\
&\qquad+\Big.x_{-1}^{-2}\{1+x_{-2}^{-2}\{1+\dots x_{-L}^{-2}\{1+x_{-L-1}^{-2}\}\dots\}\}\Big)\\
\end{align*}

as desired. 
The following two relations are straightforward computations:

\begin{equation*}
x_1^{-2}\{1+x_2^{-2}\{1+\dots x_{L-1}^{-2}\{1+x_L^{-2}\}\dots\}\}=\sum_{n=1}^{L}\frac{\varphi_{\bar b}^{L,k}(n)^2}{\varphi_{\bar b}^{L, k}(0)^2}
\end{equation*}

\begin{equation*}
x_{-1}^{-2}\{1+x_{-2}^{-2}\{1+\dots x_{-L+1}^{-2}\{1+x_{-L}^{-2}\}\dots\}\}=\sum_{n=-1}^{-L}\frac{\varphi_{\bar b}^{L,k}(n)^2}{\varphi_{\bar b}^{L, k}(0)^2}
\end{equation*}

Thus, using the fact that the eigenfunctions are normalized, we get the second expression for the determinant in $(\ref{eq3})$.

We also note that $$\left|x_1^{-1}\dots x_m^{-1}\right|=\left|\varphi_{\bar b}^{L,k}(0)\right|^{-1} \left|\varphi_{\bar b}^{L,k}(m)\right|.$$

Now, we are in a position to carry out the substitution:

\begin{align*}\rho_L(m,0)&=\int\dots\int\left(\sum_{k}\left|\langle{\delta_m,\varphi_{\bar{b}}^{L,k}}\rangle\right|\left|\langle{ \delta_0,\varphi_{\bar b}^{L,k}}\rangle\right|\right)\prod_{n=-L}^Lr_n(b_n-c_n)db_{-L}\dots db_L\\
&=\sum_{k}\int\dots\int\left|\varphi_{\bar{b}}^{L,k}(m)\right|\left|\varphi_{\bar b}^{L,k}(0)\right|\prod_{n=-L}^Lr_n(b_n-c_n)db_{-L}\dots db_L\\
&=\left(\prod_{n=-L}^{L-1}a_n\right)\sum_{k}\int\dots\int\left|\varphi_{\bar{b}}^{L,k}(m)\right|\left|\varphi_{\bar b}^{L,k}(0)\right|^{-1}\prod_{n=-L}^Lr_n(b_n-c_n)\left(\prod_{n=-L}^{L-1}a_n\right)^{-1}\left|\varphi_{\bar b}^{L,k}(0)\right|^2db_{-L}\dots db_L\\
&\leq \left(\prod_{n=-L}^{L-1}a_n\right)\int_{\Sigma_0}\int_{\mathbb{R}^{2L}}\left|x_1^{-1}\dots x_m^{-1}\right|\left(\prod_{n=-1}^{-L}r_n(E-a_{n-1}x_{n-1}^{-1}-a_nx_n-c_n)\right)r_n(E-a_{-1}x_{-1}^{-1}-a_0x_1^{-1}-c_0)\\
&\hspace{3cm}\times \left(\prod_{n=1}^{L}r_n(E-a_{n}x_{n+1}^{-1}-a_{n-1}x_n-c_n)\right)dx_{-L}\dots dx_{-1}dx_1\dots dx_LdE
\end{align*}

Let $\phi_{k;E;m}^{(n)}(x)=r_k(E-c_m-a_nx)$. Then, a quick computation shows: 

\begin{align*}\left(S_{E;0}^{(0)}\dots S_{E;-L+1}^{(-L+1)}\phi_{-L;E;-L}^{(-L)}\right)(x_1)&=\left(\prod_{n=0}^{-L+1}a_n\right)\int_{\mathbb{R}^L}r_0(E-a_{-1}x_{-1}^{-1}-a_0x_1-c_0)\\
&\hspace{2cm}\times \prod_{n=-1}^{-L}r_n(E-a_{n-1}x_{n-1}^{-1}-a_nx_n-c_n)dx_{-1}\dots dx_{-L}.
\end{align*}
Thus,

\begin{align*}
\left(US_{E;0}^{(0)}\dots S_{E;-L+1}^{(-L+1)}\phi_{-L;E;-L}^{(-L)}\right)(x_1)&=\left(\prod_{n=0}^{-L+1}a_n\right)\int_{\mathbb{R}^L}|x_1|^{-1}r_0(E-a_{-1}x_{-1}^{-1}-a_0x_1^{-1}-c_0)\\
&\hspace{2cm}\times \prod_{n=-1}^{-L}r_n(E-a_{n-1}x_{n-1}^{-1}-a_nx_n-c_n)dx_{-1}\dots dx_{-L}.
\end{align*}

Similarly,

\begin{align*}
\left(T_{E;1}^{(1)}\dots T_{E;m-1}^{(m-1)}S_{E;m}^{(m)}\dots S_{E;L-1}^{(L-1)}\phi_{L;E;L}^{(L-1)}\right)(x_1)&=
\frac{\sqrt{a_0a_{m-1}}}{a_{L-1}}\left(\prod_{n=1}^{L-1}a_n\right)\int_{\mathbb{R}^{L-1}}\left|x_2^{-1}\dots x_m^{-1}\right|\\
&\hspace{3cm}\times \prod_{n=1}^{L}r_n(E-a_{n}x_{n+1}^{-1}-a_{n-1}x_n-c_n)dx_{L}\dots dx_{2}.
\end{align*}

Combining these results, we have thus proved the following lemma:

\begin{lemm}With notation as above we have 
\begin{equation*}
\rho_L(m,0)\leq\frac{\sqrt{a_0a_{m-1}}}{a_{-L}a_{L-1}}\int_{\Sigma_0}\Big\langle{T_{E;1}^{(1)}\dots T_{E;m-1}^{(m-1)}S_{E;m}^{(m)}\dots S_{E;L-1}^{(L-1)}\phi_{L;E;L}^{(L-1)}, US_{E;0}^{(0)}\dots S_{E;-L+1}^{(-L+1)}\phi_{-L;E;-L}^{(-L)}}\Big\rangle_{L^2(\mathbb{R}, dx_1)}dE.
\end{equation*}
\end{lemm}
\section{Norm Estimates}

\begin{defi} The norm of an operator $A:L^p(\mathbb{R})\to L^q(\mathbb{R})$ will be denoted by $\norm{A}_{p,q}$.
\end{defi}
\begin{rmk}We want to point out that the following results hold for any $\alpha\in \mathbb{R}$, but since we will eventually care only for $\alpha\in \Sigma_0$ we state them in this form.
\end{rmk}
\begin{lemm}For all $\alpha \in \Sigma_0$, we have  $$\norm{S_{\alpha}^{(n)}}_{1,1}\leq 1,$$ for all $n$.
\end{lemm}

\begin{proof} We prove the statement for $n>0$,  the cases $n=0$ and $n<0$ are proved similarly. For $f\in L^1(\mathbb{R})$ we have:

\begin{align*}
\norm{S_{\alpha}^{(n)}f}_1&=\int \left|\left(S_{\alpha}^{(n)}f\right)(x)\right|dx\\
&=\int \left|a_{n-1}\int r_n(\alpha-a_{n-1}x-a_ny^{-1})f(y)dy\right|dx\\
&\leq a_{n-1}\int\int\left|d_n^{-1}r\left(d_n^{-1}\left(\alpha-a_{n-1}x-a_ny^{-1}\right)\right)\right||f(y)|dydx\\
&=\frac{a_{n-1}}{d_{n}}\int\left(\int r\left(d_n^{-1}\left(\alpha-a_{n-1}x-a_ny^{-1}\right)\right)dx\right)|f(y)|dy\\
&=\frac{a_{n-1}}{d_{n}}\int\left(\frac{d_n}{a_{n-1}}\int r(\bar x)d\bar x\right)|f(y)|dy\\
&=\norm{f}_1.
\end{align*}
We have used the fact that $r$ is nonnegative and $\norm{r}_1=1.$
\end{proof}

\begin{lemm}
For all $\alpha \in \Sigma_0$ and all $n$ we have $$\norm{S_{\alpha}^{(n)}}_{1,2}\leq \sqrt{d_n^{-1}a_{n-1}\norm{r}_{\infty}}<\infty$$
\end{lemm}
\begin{proof} We prove for the case $n>0$, the cases $n=0$ and $n<0$ are proved similarly.
For $f\in L^1(\mathbb{R}),$ we have
\begin{align*}
\norm{S_\alpha^{(n)}f}_2^2&=\int\left|\left(S_\alpha^{(n)}f\right)\right|^2dx\\
&=\int\left|\left(a_{n-1}\int r_n(\alpha -a_{n-1}x-a_ny^{-1})f(y)dy\right)\left(a_{n-1}\int r_n(\alpha-a_{n-1}x-a_nz^{-1})f(z)dz\right)\right|dx\\
&\leq \frac{a_{n-1}^2}{d_n}\norm{r}_\infty\int\left(\int|f(y)|dy\right)\left(\int r_n(\alpha-a_{n-1}x-a_nz^{-1})|f(z)|dz\right)dx\\
&=\frac{a_{n-1}^2}{d_n}\norm{r}_\infty\norm{f}_1\int\int r_n(\alpha-a_{n-1}x-a_nz^{-1})|f(z)|dzdx\\
&=\frac{a_{n-1}^2}{d_n}\norm{r}_\infty\norm{f}_1\int\left(\int r_n(\alpha-a_{n-1}x-a_nz^{-1})dx\right)|f(z)|dz\\
&=\frac{a_{n-1}^2}{d_n}\norm{r}_\infty\norm{f}_1\int\frac{1}{a_{n-1}}\left(\int r_n(\bar x)d\bar x\right)|f(z)|dz\\
&=\frac{a_{n-1}}{d_n}\norm{r}_\infty\norm{f}_1^2.
\end{align*}

So, $$\norm{S_\alpha^{(n)}f}_2\leq\sqrt{d_n^{-1}a_{n-1}\norm{r}_\infty}\,\norm{f}_1.$$
\end{proof}

\begin{lemm}\label{lem3.5}
For all $\alpha\in \Sigma_0$ we have $$\norm{T_\alpha^{(n)}}_{2,2}\leq 1.$$
\end{lemm}

\begin{proof}
Define an operator $\bar U^{(n)}$ by $$\left(\bar U^{(n)}f\right)(x)=\sqrt{\frac{a_n}{a_{n-1}}}|x|^{-1}f\left(-\frac{a_n}{a_{n-1}}x^{-1}\right).$$ We first note that $\bar U^{(n)}$ is a unitary operator on $L^2(\mathbb{R}).$ Indeed, for any $f\in L^2(\mathbb{R})$, we have

\begin{align*}
\norm{\bar U^{(n)}f}_2^2&=\int \left|\left(\bar U^{(n)}f\right)(x)\right|^2dx\\
&=\int\left|\sqrt{\frac{a_n}{a_{n-1}}}|x|^{-1}f\left(-\frac{a_n}{a_{n-1}}x^{-1}\right)\right|^2dx\\
&=\frac{a_n}{a_{n-1}}\int\left|\frac{a_{n-1}}{a_n}|u|f(u)\right|^2\frac{a_n}{a_{n-1}}u^{-2}du\\
&=\int|f(u)|^2du\\
&=\norm{f}_2^2.
\end{align*}
In the second line we have used the substitution $\displaystyle u=-\frac{a_n}{a_{n-1}}x^{-1}.$
Next, let us define an operator $K_{k;\alpha}^{(n)}$ by $\displaystyle K_{k;\alpha}^{(n)}f=r_{k;\alpha}^{(n)}\ast f;$ that is
\begin{align*}\left(K_{k;\alpha}^{(n)}f\right)(x)&=\left(r_{k;\alpha}^{(n)}\ast f\right)(x)\\
&=\int r_{k;\alpha}^{(n)}(x-y)f(y)dy\\
&=\int r_k(\alpha-a_{n-1}x+a_{n-1}y)f(y)dy.
\end{align*}

Then,

\begin{align*}
\left(K_{n;\alpha}^{(n)}\bar U^{(n)}f\right)(x)&=\left(r_{n;\alpha}^{(n)}\ast \bar U^{(n)}f\right)(x)\\
&=\int r_{n;\alpha}^{(n)}(x-y)\left(\bar U^{(n)}f\right)(y)dy\\
&=\int r_n(\alpha-a_{n-1}x+a_{n-1}y) \sqrt{\frac{a_n}{a_{n-1}}}|y|^{-1}f\left(-\frac{a_n}{a_{n-1}}y^{-1}\right)dy\\
&= \sqrt{\frac{a_n}{a_{n-1}}}\int r_n(\alpha-a_{n-1}x-a_nu^{-1})\frac{a_{n-1}}{a_n}|u|f\left(u\right)\frac{a_n}{a_{n-1}}u^{-2}du\\
&= \sqrt{\frac{a_n}{a_{n-1}}}\int r_n(\alpha-a_{n-1}x-a_nu^{-1})|u|^{-1}f\left(u\right)du\\
&=\frac{1}{a_{n-1}}\left(T_\alpha^{(n)}f\right)(x).\\
\end{align*}
We have used the substitution $\displaystyle u=-\frac{a_n}{a_{n-1}}y^{-1}.$
So, we have $\displaystyle T_\alpha^{(n)}=a_{n-1}K_{n;\alpha}^{(n)}\bar U^{(n)}.$ Then, it follows

\begin{align*}
\norm{T_\alpha^{(n)}f}_2&=\norm{a_{n-1}K_{n;\alpha}^{(n)}\bar U^{(n)}f}_2\\
&=a_{n-1}\norm{r_{n;\alpha}^{(n)}\ast \bar U^{(n)}f}_2\\
&=a_{n-1}\norm{\widehat{r_{n;\alpha}^{(n)}\ast \bar U^{(n)}f}}_2\\
&=a_{n-1}\norm{\widehat{r_{n;\alpha}^{(n)}}\cdot\widehat{ \bar U^{(n)}f}}_2\\
&\leq a_{n-1}\norm{\widehat{r_{n;\alpha}^{(n)}}}_\infty\norm{\widehat{ \bar U^{(n)}f}}_2\\
&\leq a_{n-1}\norm{r_{n;\alpha}^{(n)}}_1\norm{f}_2\\
&=a_{n-1}\norm{f}_2\int\left|r_{n;\alpha}^{(n)}(x)\right|dx\\
&=a_{n-1}\norm{f}_2\int r_n(\alpha-a_{n-1}x)dx\\
&=\norm{f}_2\int r_n(\bar x)d\bar x\\
&=\norm{r_n}_1\norm{f}_2\\
&=\norm{f}_2.
\end{align*}

Hence, the result.

\end{proof}

\begin{lemm}\label{lem3.6}
For all $\alpha,\beta \in \Sigma_0$ the operator $T_\alpha^{(n)}T_\beta^{(n+1)}$ is compact.
\end{lemm}

\begin{proof}

Let $K_{k;\alpha}^{(n)}$ and $\bar U^{(n)}$ be as before, and let $F$ be the Fourier transform, $F:f\mapsto \widehat{f}$; that is $$F[f](s)=\widehat f(s)=\int_{\mathbb{R}}e^{-2\pi i sx}f(x)dx.$$ Consider the operators $\displaystyle \bar K_{k;\alpha}^{(n)}=FK_{k;\alpha}^{(n)}F^{-1}$ and $\mathcal{U}^{(n)}=F\bar U^{(n)}F^{-1}.$
Then 
\begin{align*}
T_\alpha^{(n)}T_\beta^{(n+1)}&=\left(a_{n-1}K_{n;\alpha}^{(n)}\bar U^{(n)}\right)\left(a_{n}K_{n+1;\beta}^{(n+1)}\bar U^{(n+1)}\right)\\
&=a_{n-1}a_nF^{-1}\bar K_{n;\alpha}^{(n)}\mathcal{U}^{(n)}\bar K_{n+1;\beta}^{(n+1)}\mathcal{U}^{(n+1)}F.
\end{align*}

Since $F$ and $\mathcal{U}^{(m)}$ are unitary operators, it suffices to show that 
$\displaystyle \bar K_{n;\alpha}^{(n)}\mathcal{U}^{(n)}\bar K_{n+1;\beta}^{(n+1)}$ is compact. We will actually show that it is a Hilbert-Schmidt operator, by showing that it is an integral operator with an $L^2$ kernel, and thus compact. 
Observe that $$\bar K_{n;\alpha}^{(n)}f=FK_{n;\alpha}^{(n)}F^{-1}f=\widehat{r_{n;\alpha}^{(n)}\ast \check{ f}}=\widehat{r_{n;\alpha}^{(n)}}\cdot f.$$

Now, let $g_1\in C_c^\infty(\mathbb{R})$ be such that it is identically $1$ in some neighborhood of zero, and put $g_2=1-g_1.$ We define the following two operators

\begin{eqnarray*}
\left(U_1^{(n)}f\right)(x)&=&g_1\left(\frac{a_{n-1}}{a_n}x\right)\left(\bar U^{(n)}f\right)(x)\\
\left(U_2^{(n)}f\right)(x)&=&g_2\left(\frac{a_{n-1}}{a_n}x\right)\left(\bar U^{(n)}f\right)(x).\\
\end{eqnarray*}
Note that $\displaystyle \bar U^{(n)}=U_1^{(n)}+U_2^{(n)}.$ Then,

\begin{align*}
\left(\widehat{U_1^{(n)}f}\right)(k)&=\int e^{-2\pi ikx}\left(U_1^{(n)}f\right)(x)dx\\
&=\int e^{-2\pi ikx}g_1\left(\frac{a_{n-1}}{a_n}x\right)\left(\bar U^{(n)}f\right)(x)dx\\
&=\int e^{-2\pi ikx}g_1\left(\frac{a_{n-1}}{a_n}x\right)\sqrt{\frac{a_n}{a_{n-1}}}|x|^{-1}f\left(-\frac{a_n}{a_{n-1}}x^{-1}\right)dx\\
&=\sqrt{\frac{a_n}{a_{n-1}}}\int e^{-2\pi i\frac{a_n}{a_{n-1}}k\bar x^{-1}}g_1\left(\bar x^{-1}\right)\frac{a_{n-1}}{a_n}|\bar x|f\left(-\bar x \right)\frac{a_n}{a_{n-1}}\bar x^{-2}d\bar x\\
&=\sqrt{\frac{a_n}{a_{n-1}}}\int e^{-2\pi i\frac{a_n}{a_{n-1}}k\bar x^{-1}}g_1\left(\bar x^{-1}\right)|\bar x|^{-1}\left(\int e^{-2\pi i\bar xp} \widehat{f}(p)dp\right)d\bar x\\
&=\sqrt{\frac{a_n}{a_{n-1}}}\int\left( \int e^{-2\pi i\frac{a_n}{a_{n-1}}k\bar x^{-1}-2\pi i\bar x p}g_1\left(\bar x^{-1}\right)|\bar x|^{-1}d\bar x\right)\widehat{f}(p)dp\\
&=\sqrt{\frac{a_n}{a_{n-1}}}\int\left( \int e^{-2\pi i\frac{a_n}{a_{n-1}}kx-2\pi ipx^{-1}}g_1\left(x \right)| x|^{-1}dx\right)\widehat{f}(p)dp\\
&=\sqrt{\frac{a_n}{a_{n-1}}}\int a_1^{(n)}(k,p)\widehat{f}(p)dp,
\end{align*}
where $$a_1^{(n)}(k,p) \defeq \int e^{-2\pi i\frac{a_n}{a_{n-1}}kx-2\pi ipx^{-1}}g_1\left(x \right)| x|^{-1}dx$$

We have used the following two substitutions in this order $\displaystyle \bar x=\frac{a_n}{a_{n-1}}x^{-1}$ and $\displaystyle x=\bar x^{-1}$, in lines four and seven, respectively. 

Similarly $$\left(\widehat{U_2^{(n)}f}\right)(k)=\sqrt{\frac{a_n}{a_{n-1}}}\int a_2^{(n)}(k,p)\widehat{f}(p)dp,$$

where $$a_2^{(n)}(k,p) \defeq \int e^{-2\pi i\frac{a_n}{a_{n-1}}kx-2\pi ipx^{-1}}g_2\left(x \right)| x|^{-1}dx$$

We claim that
 \begin{equation}\label{eq5}\left (\bar K_{n;\alpha}^{(n)}\mathcal{U}^{(n)}\bar K_{n+1;\beta}^{(n+1)}f\right)(k)=\sqrt{\frac{a_n}{a_{n-1}}}\int b^{(n)}(k,p)f(p)dp,\end{equation}

where

 $$b^{(n)}(k,p)=\widehat{r_{n;\alpha}^{(n)}}(k)\Big(a_1^{(n)}(k,p)+a_2^{(n)}(k,p)\Big)\widehat{r_{n+1;\beta}^{(n+1)}}(p).$$
Observe that 
$$\bar K_{n;\alpha}^{(n)}\mathcal{U}^{(n)}\bar K_{n+1;\beta}^{(n+1)}=\bar K_{n;\alpha}^{(n)}FU_1^{(n)}F^{-1}\bar K_{n+1;\beta}^{(n+1)}+\bar K_{n;\alpha}^{(n)}FU_2^{(n)}F^{-1}\bar K_{n+1;\beta}^{(n+1)},$$

where we have used the fact that $\bar U^{(n)}=U_1^{(n)}+U_2^{(n)}.$

Next

\begin{align*}
\left(\bar K_{n;\alpha}^{(n)}FU_1^{(n)}F^{-1}\bar K_{n+1;\beta}^{(n+1)}f\right)(k)&=\widehat{r_{n;\alpha}^{(n)}}(k)\cdot\left(FU_1^{(n)}F^{-1}\bar K_{n+1;\beta}^{(n+1)}f\right)(k)\\
&=\sqrt{\frac{a_n}{a_{n-1}}}\int\widehat{r_{n;\alpha}^{(n)}}(k)a_1^{(n)}(k,p)\widehat{r_{n+1;\beta}^{(n+1)}}(p)f(p)dp.\\
\end{align*}
Similarly,

\begin{align*}
\left(\bar K_{n;\alpha}^{(n)}FU_2^{(n)}F^{-1}\bar K_{n+1;\beta}^{(n+1)}f\right)(k)&=\widehat{r_{n;\alpha}^{(n)}}(k)\cdot\left(FU_2^{(n)}F^{-1}\bar K_{n+1;\beta}^{(n+1)}f\right)(k)\\
&=\sqrt{\frac{a_n}{a_{n-1}}}\int\widehat{r_{n;\alpha}^{(n)}}(k)a_2^{(n)}(k,p)\widehat{r_{n+1;\beta}^{(n+1)}}(p)f(p)dp.\\
\end{align*}
Combining these two expressions, we get $(\ref{eq5})$. Next, we need to show that $b^{(n)}$ is in $L^2.$ We have, 

\begin{align}\label{ineq1}\nonumber\norm{b^{(n)}(k,p)}_{L^2(\mathbb{R},dk)\times L^2(\mathbb{R}, dp)}&\leq \norm{\widehat{r_{n;\alpha}^{(n)}}(k)a_1^{(n)}(k,p)\widehat{r_{n+1;\beta}^{(n+1)}}(p)}_{L^2(\mathbb{R},dk)\times L^2(\mathbb{R}, dp)}\\ 
&+\norm{\widehat{r_{n;\alpha}^{(n)}}(k)a_2^{(n)}(k,p)\widehat{r_{n+1;\beta}^{(n+1)}}(p)}_{L^2(\mathbb{R},dk)\times L^2(\mathbb{R}, dp)}.
\end{align}

Note that,

\begin{align*}
\norm{\widehat{r_{n;\alpha}^{(n)}}(k)a_1^{(n)}(k,p)\widehat{r_{n+1;\beta}^{(n+1)}}(p)}_{L^2(\mathbb{R},dk)\times L^2(\mathbb{R}, dp)}&=\left(\int_{\mathbb{R}\times\mathbb{R}}\left|\widehat{r_{n;\alpha}^{(n)}}(k)a_1^{(n)}(k,p)\widehat{r_{n+1;\beta}^{(n+1)}}(p)\right|^2dkdp\right)^{1/2}\\
&\leq \left(\int_{\mathbb{R}}\int_\mathbb{R}\left|\widehat{r_{n;\alpha}^{(n)}}(k)\right|^2\left|a_1^{(n)}(k,p)\right|^2dkdp\right)^{1/2}\norm{\widehat{r_{n+1;\beta}^{(n+1)}}}_{L^{\infty}(\mathbb{R},dp)}\\
&\leq \left(\int_{\mathbb{R}}\left|\widehat{r_{n;\alpha}^{(n)}}(k)\right|^2\left(\int_\mathbb{R}\left|a_1^{(n)}(k,p)\right|^2dp\right)dk\right)^{1/2}\norm{\widehat{r_{n+1;\beta}^{(n+1)}}}_{L^{\infty}(\mathbb{R},dp)}\\
&\leq \norm{\widehat{r_{n;\alpha}^{(n)}}}_{L^2(\mathbb{R},dk)}\sup_k\left(\int_\mathbb{R}\left|a_1^{(n)}(k,p)\right|^2dp\right)^{1/2}\norm{\widehat{r_{n+1;\beta}^{(n+1)}}}_{L^{\infty}(\mathbb{R},dp)}\\
&=\norm{\widehat{r_{n;\alpha}^{(n)}}}_{L^2(\mathbb{R},dk)}\sup_k\norm{a_1^{(n)}(k,\cdot)}_{L^2(\mathbb{R},dp)}\norm{\widehat{r_{n+1;\beta}^{(n+1)}}}_{L^{\infty}(\mathbb{R},dp)}.\\
\end{align*}

Similarly, 

\begin{align*}
\norm{\widehat{r_{n;\alpha}^{(n)}}(k)a_2^{(n)}(k,p)\widehat{r_{n+1;\beta}^{(n+1)}}(p)}_{L^2(\mathbb{R},dk)\times L^2(\mathbb{R}, dp)}&\leq \norm{\widehat{r_{n;\alpha}^{(n)}}}_{L^\infty(\mathbb{R},dk)}\sup_p\norm{a_2^{(n)}(\cdot,p)}_{L^2(\mathbb{R},dk)}\norm{\widehat{r_{n+1;\beta}^{(n+1)}}}_{L^2(\mathbb{R},dp)}.
\end{align*}

Since, $r\in L^1(\mathbb{R})\cap L^\infty(\mathbb{R})$, then $r\in L^1(\mathbb{R})\cap L^2(\mathbb{R}).$ So, by Plancherel's theorem, it follows that $\widehat{r}\in L^2(\mathbb{R})\cap L^\infty(\mathbb{R}).$ Thus, it remains to show that

\begin{enumerate}[i.]
\item $\displaystyle \sup_k\norm{a_1^{(n)}(k,\cdot)}_{L^2(\mathbb{R},dp)}<\infty$, and
\item $\displaystyle \sup_p\norm{a_2^{(n)}(\cdot,p)}_{L^2(\mathbb{R},dk)}<\infty.$
\end{enumerate}

We begin by proving the first claim. To this end let 

\begin{eqnarray*}
f_N^{(k)}(\bar x)&=&e^{-2\pi i\frac{a_n}{a_{n-1}}k\bar x^{-1}}\frac{g_1(\bar x^{-1})}{|\bar x|}\cdot \chi_{[-N,N]}(\bar x)\\
f^{(k)}(\bar x)&=&e^{-2\pi i\frac{a_n}{a_{n-1}}k\bar x^{-1}}\frac{g_1(\bar x^{-1})}{|\bar x|}.
\end{eqnarray*}

Since $g_1$ is compactly supported and is identically equal to $1$ in a neighborhood of $0$, it is not difficult to see that $f^{(k)}$ is an $L^2$ function, and that it's $L^2$ norm is independent of $k$.  From this, it is, also, not difficult to see that $f_N^{(k)}$ converges to $f^{(k)}$ in $L^2$ sense. As a result, it is straightforward to see that $\widehat{f_N^{(k)}}$ converges to $\widehat{f^{(k)}}$ in $L^2$ sense; where 
\begin{eqnarray*}
\widehat{f^{(k)}}(p)&=&\int f^{(k)}(\bar x)e^{-2\pi ip\bar x}d\bar x=\int e^{-2\pi i\frac{a_n}{a_{n-1}}k\bar x^{-1}-2\pi ip\bar x}\,\frac{g_1(\bar x^{-1})}{|\bar x|}d\bar x
\end{eqnarray*}

$$\widehat{f_N^{(k)}}=\int_{|\bar x|<N} e^{-2\pi i\frac{a_n}{a_{n-1}}k\bar x^{-1}-2\pi ip\bar x}\, \frac{g_1(\bar x^{-1})}{|\bar x|}d\bar x.$$
Note that,

\begin{align*}
a_1^{(n)}(k,p)\cdot\chi_{\{x:|x|>\frac{1}{N}\}}(x)&=\int_{|x|>\frac{1}{N}}e^{-2\pi i\frac{a_n}{a_{n-1}}kx-2\pi ipx^{-1}}g_1(x)|x|^{-1}dx\\
&=\int_{|\bar x|<N}e^{-2\pi i\frac{a_n}{a_{n-1}}k\bar x^{-1}-2\pi ip\bar x}\, \frac{g_1(\bar x^{-1})}{|\bar x|}d\bar x\\
&=\widehat{f^{(k)}_N}(p).
\end{align*}
Then, from our discussion above, it follows that $\displaystyle a_1^{(n)}(k,p)=\widehat{f^{(k)}}(p).$ Hence, by unitarity of the Fourier transform, and the fact that $f^{(k)}$ has $L^2$ norm independent of $k$, we get

\begin{equation}\label{eqqq1}\sup_k\norm{a_1^{(n)}(k,\cdot)}_{L^2(\mathbb{R},dp)}=\sup_k\norm {\widehat{f^{(k)}}}_{L^2(\mathbb{R}, dp)}=\sup_k\norm {f^{(k)}}_{L^2(\mathbb{R}, dp)}=\norm {f^{(k)}}_{L^2(\mathbb{R}, dp)}<\infty.\end{equation}

Next, let $$f^{(p)}(\bar x)=\, e^{-2\pi i\frac{a_n}{a_{n-1}}p\bar x^{-1}}\, g_2\left(\frac{a_{n-1}}{a_n}\bar x\right)|\bar x|^{-1}.$$ Since $g_2$ vanishes in a neighborhood of $0$, it is easy to see that $f^{(p)}$ is an $L^2$ function, and that it's norm is independent of $p$.  Then,

\begin{align*}
\widehat{f^{(p)}}(k)&=\int e^{-2\pi ik \bar x} f^{(p)}(\bar x)\bar x\\
&=\int e^{-2\pi ik\bar x-2\pi i\frac{a_n}{a_{n-1}}p\bar x^{-1}}\, g_2\left(\frac{a_{n-1}}{a_n}\bar x\right)|\bar x|^{-1}d\bar x\\
&=\int e^{-2\pi i\frac{a_n}{a_{n-1}}kx-2\pi ipx^{-1}}\,g_2(x)|x|^{-1}dx\\
&=a_2^{(n)}(k,p).
\end{align*}

Hence, \begin{equation}\label{eqqq2}\sup_p\norm{a_2^{(n)}(\cdot, p)}_{L^2(\mathbb{R},dk)}=\sup_p\norm{\widehat{f^{(p)}}}_{L^2(\mathbb{R}, dk)}=\sup_p\norm{f^{(p)}}_{L^2(\mathbb{R}, dk)}=\norm{f^{(p)}}_{L^2(\mathbb{R}, dk)}<\infty.\end{equation}

This concludes that $b^{(n)}$ is an $L^2$ function. So, $\displaystyle T_\alpha^{(n)}T_\beta^{(n+1)}$ is Hilbert-Schmidt, and thus compact. 
\end{proof}

Next, we adopt the technique developed in \cite{bs} to prove the following lemma.

\begin{lemm}\label{lemm4.8}
For some fixed constant $C_0$ we have $$\norm{T_\alpha^{(n)}T_\beta^{(n+1)}f}_2\leq A(n,n+1)\norm{f}_2,$$ where 

$$A(n,n+1)\defeq\left(\frac{15}{16}+\frac{1}{16}\sup_{|k|\geq t_n \frac{C_0}{\norm{a}_\infty}}\left|\widehat{r}(k)\right|^2\right)^{\frac{1}{2}},$$ where $\displaystyle t_n=\min\left(d_n, d_{n+1}\right)$.
\end{lemm}

\begin{proof}
Above we have shown that, in particular, $T_\alpha^{(n)}T_\beta^{(n+1)}$ is a Hilbert-Schmidt operator. Specifically, 

$$T_\alpha^{(n)}T_\beta^{(n+1)}=F^{-1}\left(a_{n-1}a_n\bar K_{n;\alpha}^{(n)}\mathcal{U}^{(n)}\bar K_{n+1;\beta}^{(n+1)}\right)\mathcal{U}^{(n+1)}F.$$

So, it suffices to show that for $\norm{\varphi}_2=\norm{\psi}_2=1$ we have 
\begin{equation}\label{eq10}\left|\langle{\varphi, a_{n-1}a_n\bar K_{n;\alpha}^{(n)}\mathcal{U}^{(n)}\bar K_{n+1;\beta}^{(n+1)}\psi\rangle}\right|\leq A(n,n+1).\end{equation}

Pick $C_0$, such that 
\begin{equation}\label{eq9}
B\left(\int_{|k|\leq C_0}\int_{ |p|\leq C_0}\left|a^{(n)}(k,p)\right|^2dkdp\right)^{1/2}\leq \frac{7}{16},
\end{equation}

where $B=\sup\sqrt{\frac{a_n}{a_{n-1}}}$. We claim that, this is possible, since the left hand side of $(\ref{eq9})$ goes to zero, as $C_0\to 0$, and also that such a $K_0$ can be chosen independently of $n$. Both of these facts are a byproduct of the proof of Lemma $\ref{lem3.6}.$ More precisely, note that

\begin{align}\label{eqq12}\nonumber
\norm{a^{(n)}(k,p)}_{L^2\left([-C_0, C_0]^2, dkdp\right)}&\leq\norm{a_1^{(n)}(k,p)}_{L^2\left([-C_0, C_0]^2, dkdp\right)}+\norm{a_2^{(n)}(k,p)}_{L^2\left([-C   _0, C_0]^2, dkdp\right)}\\ \nonumber
&\leq \sqrt{2C_0}\left(\sup_{k}\norm{a_1^{(n)}(k,\cdot)}_{L^2\left([-C_0, C_0], dp\right)}+\sup_{p}\norm{a_2^{(n)}(\cdot,p)}_{L^2\left([-C_0, C_0], dk\right)}\right)\\ \nonumber
&\leq\sqrt{2C_0}\left(\sup_{k}\norm{a_1^{(n)}(k,\cdot)}_{L^2\left(\mathbb{R}, dp\right)}+\sup_{p}\norm{a_2^{(n)}(\cdot,p)}_{L^2\left(\mathbb{R}, dk\right)}\right)\\ \nonumber
&=\sqrt{2C_0}\left(\norm{f^{(k)}}_{L^2(\mathbb{R},dp)}+\norm{f^{(p)}}_{L^2(\mathbb{R}, dk)}\right)\\
&\leq \sqrt{2C_0}\times \sqrt{2\pi}\left(\left(\int_{\mathbb{R}}\frac{\left|g_1(p^{-1})\right|^2}{|p|^2}dp\right)^{1/2}+\sqrt{\frac{\norm{a}_\infty}{\delta}}\left(\int_{\mathbb{R}}\frac{\left|g_2(y)\right|^2}{|y|^2}dy\right)^{1/2}\right),
\end{align}
where, going from line three to four, we have used expressions $(\ref{eqqq1})$ and $(\ref{eqqq2})$, and from line four to five we have performed a change of variables and used the fact that $0<\delta\leq a_n\leq \norm{a}_\infty,$ for all $n.$ So, using the fact that, as seen before, the integrals that appear above are finite, we can pick $K_0$ independently of $n$, such that the right hand side of $(\ref{eqq12})$ is less than $\frac{7}{16}.$

Let $\varphi_+=\varphi\chi_{\{|k|\geq C_0\}}$ and $\psi_+=\psi\chi_{\{|k|\geq C_0\}}.$ We consider two cases 

\begin{enumerate}[(i)]
\item $\norm{\varphi_+}_2\geq \frac{1}{4}$ or $\norm{\psi_+}_2\geq \frac{1}{4};$
\item $\norm{\varphi_+}_2\leq \frac{1}{4}$ and $\norm{\psi_+}_2\leq \frac{1}{4}.$
\end{enumerate}
 First, using the fact that $\bar K_{n;\alpha}^{(n)}f=\widehat{r^{(n)}_{n;\alpha}}\cdot f$ and the fact that $\|\widehat{g}\|_\infty\leq \|g\|_1,$ we get that \begin{equation}\label{eq13}\norm{\bar K_{n;\alpha}^{(n)}f}_2\leq \frac{1}{a_{n-1}}\norm{f}_2,\end{equation} and also

\begin{align*}
\widehat{r_{n;\alpha}^{(n)}}(s)&=\int e^{-2\pi isx}r_n(\alpha-a_{n-1}x)dx\\
&=\frac{1}{a_{n-1}}e^{-2\pi is\frac{\alpha}{a_{n-1}}}\int e^{-2\pi i\left(-\frac{d_n}{a_{n-1}}s\right)\bar x}\, r(\bar x)d\bar x\\
&=\frac{1}{a_{n-1}}e^{-2\pi is\frac{\alpha}{a_{n-1}}}\,\widehat{r}\left(-\frac{d_n}{a_{n-1}}s\right).
\end{align*}
So, \begin{equation}\label{eq11}\left|\widehat{r_{n;\alpha}^{(n)}}(s)\right|=\frac{1}{a_{n-1}}\left|\widehat{r}\left(-\frac{d_n}{a_{n-1}}s\right)\right|.\end{equation}

Now, suppose that $\norm{\psi_+}_2\geq \frac{1}{4}.$ Then,

\begin{align*}
\norm{a_n\bar K_{n+1;\beta}^{(n+1)}\psi}_2^2&=a_n^2\int_\mathbb{R}\left|\left(\bar K_{n+1;\beta}^{(n+1)}\psi\right)(k)\right|^2dk\\
&=a_n^2\int_\mathbb{R}\left|\widehat r_{n+1;\beta}^{(n+1)}(k)\psi(k)\right|^2dk\\
&=a_n^2\int_\mathbb{R}\left|\widehat r_{n+1;\beta}^{(n+1)}(k)\right|^2\left|\psi(k)\right|^2dk\\
&=a_n^2\int_\mathbb{R}\left|\frac{1}{a_n}\widehat r\left(-\frac{d_{n+1}}{a_n}k\right)\right|^2\left|\psi(k)\right|^2dk\\
&=\int_\mathbb{R}\left|\widehat r\left(-\frac{d_{n+1}}{a_n}k\right)\right|^2\left|\psi(k)\right|^2dk\\
&=\int_{\{|k|\geq K_0\}}\left|\widehat r\left(-\frac{d_{n+1}}{a_n}k\right)\right|^2\left|\psi(k)\right|^2dk+\int_{\{|k|<C_0\}}\left|\widehat r\left(-\frac{d_{n+1}}{a_n}k\right)\right|^2\left|\psi(k)\right|^2dk\\
&\leq \sup_{|k|\geq C_0}\left|\widehat r\left(-\frac{d_{n+1}}{a_n}k\right)\right|^2\int_{\{|k|\geq C_0\}}\left|\psi(k)\right|^2dk+\int_{\{|k|<C_0\}}\left|\psi(k)\right|^2dk\\
&=\sup_{|k|\geq \frac{d_{n+1}}{a_n}C_0}\left|\widehat r\left(k\right)\right|^2\int_{\{|k|\geq C_0\}}\left|\psi(k)\right|^2dk+\int_{\{|k|<C_0\}}\left|\psi(k)\right|^2dk\\
&\leq\sup_{|k|\geq \frac{C_0}{\norm{a}_\infty}d_{n+1}}\left|\widehat r\left(k\right)\right|^2\int_{\{|k|\geq C_0\}}\left|\psi(k)\right|^2dk+\int_{\{|k|<C_0\}}\left|\psi(k)\right|^2dk\\
&\hspace{6.4cm}+\int_{\{|k|\geq C_0\}}\left|\psi(k)\right|^2dk-\int_{\{|k|\geq C_0\}}\left|\psi(k)\right|^2dk\\
&=\left(\sup_{|k|\geq \frac{C_0}{\norm{a}_\infty}d_{n+1}}\left|\widehat r\left(k\right)\right|^2-1\right)\int_{\{|k|\geq C_0\}}\left|\psi(k)\right|^2dk+1\\
&=1+\left(-\int_{\{|k|\geq C_0\}}\left|\psi(k)\right|^2dk\right)\left(1-\sup_{|k|\geq \frac{C_0}{\norm{a}_\infty}d_{n+1}}\left|\widehat r\left(k\right)\right|^2\right)\\
&\leq 1-\frac{1}{16}\left(1-\sup_{|k|\geq \frac{C_0}{\norm{a}_\infty}d_{n+1}}\left|\widehat r\left(k\right)\right|^2\right)\\
&=1-\frac{1}{16}+\frac{1}{16}\sup_{|k|\geq \frac{C_0}{\norm{a}_\infty}d_{n+1}}\left|\widehat r\left(k\right)\right|^2\\
&=\frac{15}{16}+\frac{1}{16}\sup_{|k|\geq \frac{C_0}{\norm{a}_\infty}d_{n+1}}\left|\widehat r\left(k\right)\right|^2.\\
\end{align*}
Now, using Cauchy-Schwarz, $(\ref{eq13})$, and the fact that $\mathcal{U}^{(n)}$ is unitary, we get

\begin{align*}
\left|\langle{\varphi, a_{n-1}a_n\bar K_{n;\alpha}^{(n)}\mathcal{U}^{(n)}\bar K_{n+1;\beta}^{(n+1)}\psi\rangle}\right|&\leq\norm{\varphi}_2\norm{a_{n-1}a_n\bar K_{n;\alpha}^{(n)}\mathcal{U}^{(n)}\bar K_{n+1;\beta}^{(n+1)}\psi\rangle}_2\\
&\leq \norm{a_n\bar K_{n+1;\beta}^{(n+1)}\psi}_2.
\end{align*}
Thus, from above, in this case the result follows. 

Next, if $\norm{\varphi_+}_2\geq \frac{1}{4}$, then 

\begin{align*}
\left|\langle{\varphi, a_{n-1}a_n\bar K_{n;\alpha}^{(n)}\mathcal{U}^{(n)}\bar K_{n+1;\beta}^{(n+1)}\psi\rangle}\right|&=\left|\int_\mathbb{R} \varphi(k)a_{n-1}a_n\left(\bar K_{n;\alpha}^{(n)}\mathcal{U}^{(n)}\bar K_{n+1;\beta}^{(n+1)}\psi\right)(k)dk\right|\\
&\leq \int_\mathbb{R}\left|\varphi(k)a_{n-1}\widehat r_{n;\alpha}^{(n)}(k)\right|\left|a_n\left(\mathcal{U}^{(n)}\bar K_{n+1;\beta}^{(n+1)}\psi\right)(k)\right|dk\\
&\leq \left(\int_\mathbb{R} \left|\varphi(k)a_{n-1}\widehat r_{n;\alpha}^{(n)}(k)\right|^2dk\right)^{1/2}\left(\int_\mathbb{R} \left|a_n\left(\mathcal{U}^{(n)}\bar K_{n+1;\beta}^{(n+1)}\psi\right)(k)\right|^2dk\right)^{1/2}\\
&=\left(\int_\mathbb{R} \left|\widehat r\left(-\frac{d_n}{a_{n-1}}k\right)\right|^2|\varphi(k)|^2dk\right)^{1/2}a_n\norm{\mathcal{U}^{(n)}\bar K_{n+1;\beta}^{(n+1)}\psi}_2\\
&=\left(\int_\mathbb{R} \left|\widehat r\left(-\frac{d_n}{a_{n-1}}k\right)\right|^2|\varphi(k)|^2dk\right)^{1/2}a_n\norm{\bar K_{n+1;\beta}^{(n+1)}\psi}_2\\
&\leq \left(\int_\mathbb{R} \left|\widehat r\left(-\frac{d_n}{a_{n-1}}k\right)\right|^2|\varphi(k)|^2dk\right)^{1/2}\\
&\leq \left(\frac{15}{16}+\frac{1}{16}\sup_{|k|\geq \frac{C_0}{\norm{a}_\infty}d_n}\left|\widehat r\left(k\right)\right|^2\right)^{1/2}.
\end{align*}
The last inequality follows via the same argument as before. Thus, again, the result follows. 

\noindent Before we consider the second case, let $\varphi_-=\varphi\chi_{\{|k|<C_0\}},$ and $\psi_-=\psi\chi_{\{|k|<C_0\}}.$ Then

\begin{align*}
\left|\langle{\varphi_-, a_{n-1}a_n\bar K_{n;\alpha}^{(n)}\mathcal{U}^{(n)}\bar K_{n+1;\beta}^{(n+1)}\psi_-\rangle}\right|&=\left|\int_\mathbb{R} \varphi_-(k)a_{n-1}a_n\left(\bar K_{n;\alpha}^{(n)}\mathcal{U}^{(n)}\bar K_{n+1;\beta}^{(n+1)}\psi_-\right)(k)dk\right|\\
&=\left|\int_\mathbb{R}\varphi_-(k)\sqrt{\frac{a_n}{a_{n-1}}}a_{n-1}a_n\int_\mathbb{R}
\widehat r_{n;\alpha}^{(n)}(k)a^{(n)}(k,p)\widehat r_{n+1;\beta}^{(n+1)}(p)\psi_-(p)dpdk\right|\\
&=\left|\int_{\{|k|\leq C_0\}}\int_{\{|p|\leq C_0\}}\sqrt{\frac{a_n}{a_{n-1}}}a_{n-1}a_n\varphi(k)
\widehat r_{n;\alpha}^{(n)}(k)a^{(n)}(k,p)\widehat r_{n+1;\beta}^{(n+1)}(p)\psi(p)dpdk\right|\\
&\leq \int_{\{|k|\leq C_0\}}\int_{\{|p|\leq C_0\}}\sqrt{\frac{a_n}{a_{n-1}}}a_{n-1}a_n|\varphi(k)\psi(p)|
\left|\widehat r_{n;\alpha}^{(n)}(k)\right|\left|a^{(n)}(k,p)\right|\left|\widehat r_{n+1;\beta}^{(n+1)}(p)\right|dpdk\\
&\leq\sqrt{\frac{a_n}{a_{n-1}}} \int_{\{|k|\leq C_0\}}\int_{\{|p|\leq C_0\}}|\varphi(k)\psi(p)|
\left|a^{(n)}(k,p)\right|dpdk\\
&\leq\sqrt{\frac{a_n}{a_{n-1}}}\left( \int_{\{|k|\leq C_0\}}\int_{\{|p|\leq C_0\}}|\varphi(k)|^2|\psi(p)|^2
dpdk\right)^{1/2}\\
&\hspace{4cm}\times\left( \int_{\{|k|\leq C_0\}}\int_{\{|p|\leq C_0\}}
\left|a^{(n)}(k,p)\right|^2dpdk\right)^{1/2}\\
&=\sqrt{\frac{a_n}{a_{n-1}}}\left( \int_{\{|p|\leq C_0\}}|\psi(p)|^2
dp\right)^{1/2}\left( \int_{\{|k|\leq C_0\}}|\varphi(k)|^2
dk\right)^{1/2}\\
&\hspace{4cm}\times\left( \int_{\{|k|\leq C_0\}}\int_{\{|p|\leq C_0\}}
\left|a^{(n)}(k,p)\right|^2dpdk\right)^{1/2}\\
&\leq B\left( \int_{\{|k|\leq C_0\}}\int_{\{|p|\leq C_0\}}
\left|a^{(n)}(k,p)\right|^2dpdk\right)^{1/2}\\
&\leq\frac{7}{16}.
\end{align*}

Finally, if $\norm{\varphi_+}_2\leq \frac{1}{4}$, and $\norm{\psi_+}_2\leq \frac{1}{4}$, we have

\begin{align*}
\left|\langle{\varphi, a_{n-1}a_n\bar K_{n;\alpha}^{(n)}\mathcal{U}^{(n)}\bar K_{n+1;\beta}^{(n+1)}\psi\rangle}\right|
&=\left|\langle{\varphi_+, a_{n-1}a_n\bar K_{n;\alpha}^{(n)}\mathcal{U}^{(n)}\bar K_{n+1;\beta}^{(n+1)}\psi\rangle}+\langle{\varphi_-, a_{n-1}a_n\bar K_{n;\alpha}^{(n)}\mathcal{U}^{(n)}\bar K_{n+1;\beta}^{(n+1)}\psi\rangle}\right|\\
&\leq \norm{\varphi_+}_2+\left|\langle{\varphi_-, a_{n-1}a_n\bar K_{n;\alpha}^{(n)}\mathcal{U}^{(n)}\bar K_{n+1;\beta}^{(n+1)}\psi_+\rangle}\right|+\left|\langle{\varphi_-, a_{n-1}a_n\bar K_\alpha^{(n)}\mathcal{U}^{(n)}\bar K_\beta^{(n+1)}\psi_-\rangle}\right|\\
&\leq \norm{\varphi_+}_2+\norm{\psi_+}_2+\left|\langle{\varphi_-, a_{n-1}a_n\bar K_{n;\alpha}^{(n)}\mathcal{U}^{(n)}\bar K_{n+1;\beta}^{(n+1)}\psi_-\rangle}\right|\\
&\leq \frac{1}{4}+\frac{1}{4}+\frac{7}{16}\\
&=\frac{15}{16}< A(n,n+1).
\end{align*}
This concludes the proof of the lemma. 

\end{proof}

We record the following as a corollary, so we can refer to it later. 
\begin{cor} Let $d_n$ be such that $d_n=1$ for all $n$. Then, there exists some constant $0<q<1$, such that $$\norm{T_\alpha^{(n)}T_\beta^{(n+1)}}_{2,2}\leq q$$ for all $n$ and all $\alpha,\beta.$
\end{cor}

\begin{proof}
This is an immediate consequence of Lemma $\ref{lemm4.8}$ with $$q\defeq\left(\frac{15}{16}+\frac{1}{16}\sup_{|k|\geq \frac{C_0}{\norm{a}_\infty}}\left|\widehat{r}(k)\right|^2\right)^{\frac{1}{2}}.$$

\end{proof}

If the sequence $d_n\not\equiv 1$ we can no longer bound $\norm{T_\alpha^{(n)}T_\beta^{(n+1)}}_{2,2}$ uniformly away from $1$, however, we can still control the rate at which this norm converges to $1$, as is established in the following Lemma.  
\begin{lemm}\label{lemm4.9}
 Let  $d_n$ be a fixed sequence with $0\leq d_n\leq1$ and $d_n\geq C|n|^{-\zeta}$ for $\zeta<\frac{1}{2}$, and some constant $C>0$. Then,
$$A(s)\defeq\left(\frac{15}{16}+\frac{1}{16}\sup_{|k|\geq t_s \frac{C_0}{\norm{a}_\infty}}\left|\widehat{r}(k)\right|^2\right)^{\frac{1}{2}}\leq \exp\left(-\gamma' |s|^{-2\zeta}\right),$$ for some $\gamma'>0$, where $\displaystyle t_s=\min\left(d_{2s-1}, d_{2s}\right)$.

\end{lemm}

\begin{proof}
First let us show that $\frac{d^2}{dk^2}\left|\widehat r(k)\right|^2\Big|_{k=0}<0.$ We compute,

\begin{align*}
\frac{d^2}{dk^2}\left|\widehat r(k)\right|^2\Big|_{k=0}&=\frac{d^2}{dk^2}\widehat r(k)\overline{\widehat r(k)}\Big|_{k=0}\\
&=\frac{d}{dk}\left(\widehat r(k)\frac{d}{dk}\overline{\widehat r(k)}+\overline{\widehat r(k)}\frac{d}{dk}\widehat r(k)\right)\Big|_{k=0}\\
&=2\frac{d}{dk}\widehat r(k)\frac{d}{dk}\overline{\widehat r(k)}\Big|_{k=0}+\widehat r(k)\frac{d^2}{dk^2}\overline{\widehat r(k)}\Big|_{k=0}+\overline{\widehat r(k)}\frac{d^2}{dk^2}\widehat r(k)\Big|_{k=0}\\
&=8\pi^2\int xe^{-2\pi ikx}r(x)dx\int x e^{2\pi ikx}r(x)dx\Big|_{k=0}-8\pi^2\, \Re\left(\int e^{-2\pi ikx}r(x)dx\int x^2e^{2\pi ikx}r(x)dx\right)\Big|_{k=0}\\
&=8\pi^2 \left(\int xr(x)dx\right)^2-8\pi^2 \left(\int x^2r(x)dx\right)\\
&<0,
\end{align*}

where the strict inequality, in the last line, follows by Cauchy-Schwarz. 
Now, before we proceed, let us prove the following claim.

\begin{claim}
For $k\neq 0$ we have $$\left|\widehat r(k)\right|<1.$$
\end{claim} 

\begin{proof}[Proof of Claim.]
First, we know that in general we have $|\widehat r(k)|\leq \|r\|_1=1$ with $|\widehat r(0)|=1.$ So, suppose that there is some $k\neq 0$ such that $\left|\widehat r(k)\right|=1.$ Then, there is some $\theta(k)\in [0,2\pi)$ such that $$\widehat r(k)=e^{i\theta(k)}.$$ Then, in particular it follows that  $$\int \cos(2\pi kx)r(x)dx=\cos \theta(k).$$ Equivalently, $$\int \left(\cos(2\pi kx)-\cos\theta(k)\right)r(x)dx=0,$$ from which it follows that $$\cos(2\pi kx)=\cos\theta(k),$$ for $Leb-a.e.\, x$, which is clearly not possible.  
\end{proof}

Thus, since $\left|\widehat r(k)\right|<1$ for $k\neq 0$ and $\left|\widehat r(0)\right|=1$, with $\frac{d^2}{dk^2}\left|\widehat r(k)\right|^2\Big|_{k=0}<0,$ by a Taylor series expansion around zero, for $\lambda$ small enough we have 
$$\sup_{|k|\geq \lambda}\left|\widehat r(k)\right|^2\leq 1-c\lambda^2\leq e^{-\tilde c\lambda^2},$$ with $c\defeq-\frac{d^2}{dk^2}\left|\widehat r(k)\right|^2\Big|_{k=0}>0,$ and some $\tilde c>0.$ Then,

$$\left(\frac{15}{16}+\frac{1}{16}\sup_{|k|\geq \lambda}\left|\widehat r(k)\right|^2\right)^{1/2}\leq\left( \frac{15}{16}+\frac{1}{16}-c_1\lambda^2\right)^{1/2}=\left(1-c_1\lambda^2\right)^{1/2}\leq e^{-c_2\lambda^2}.$$ As a result, up to possibly shrinking $C_0$, we have

\begin{align*}
A(s)&\defeq\left(\frac{15}{16}+\frac{1}{16}\sup_{|k|\geq t_s \frac{C_0}{\norm{a}_\infty}}\left|\widehat{r}(k)\right|^2\right)^{\frac{1}{2}}\\
&\leq \exp\left(-c_2\left(t_s\frac{C_0}{\norm{a}_\infty}\right)^2\right)\\
&\leq \exp\left(-\gamma'|s|^{-2\zeta}\right).
\end{align*}

\end{proof}

The following corollary is an immediate consequence of the arguments above.
\begin{cor}
 Let  $d_n$ be a fixed sequence with $0\leq d_n\leq1$ and $d_n\geq C|n|^{-\zeta}$ for $\zeta<\frac{1}{2}$. Then, 
$$\norm{T_\alpha^{(2s-1)}T_\beta^{(2s)}}_{2,2}\leq \exp\left(-\gamma'|s|^{-2\zeta}\right).$$
\end{cor}

\begin{proof}
From Lemma \ref{lemm4.8} we clearly have, $$\norm{T_\alpha^{(2s-1)}T_\beta^{(2s)}}_{2,2}\leq A(s),$$ hence the result follows from Lemma \ref{lemm4.9}. 

\end{proof}

\begin{lemm}With notation as above we have $$A(1)\times \dots \times A(s)\leq \exp\left(-\gamma' |s|^{1-2\zeta}\right),$$ for some constant $\gamma'>0.$
\end{lemm}

\begin{proof}
We have,

\begin{align*}A(1)\times \dots \times A(s)&=\exp\left(-\gamma'\left(1+2^{-2\zeta}+\dots+(s-1)^{-2\zeta}+s^{-2\zeta}\right)\right)\\
&=\exp\left(-\gamma' s^{-2\zeta}\left(s^{2\zeta}+2^{-2\zeta}s^{2\zeta}+\dots+(s-1)^{-2\zeta}s^{2\zeta}+1\right)\right)\\
&\leq\exp\left(- \gamma' s^{1-2\zeta}\right).
\end{align*}
The last inequality follows from the fact that $$s^{2\zeta}+2^{-2\zeta}s^{2\zeta}+\dots+(s-1)^{-2\zeta}s^{2\zeta}+1\geq 1+2^{-2\zeta}2^{2\zeta}+\dots+(s-1)^{-2\zeta}(s-1)^{2\zeta}+1=s.$$

\end{proof}



\subsection{Proof of Theorem \ref{thm2.1}.}

\begin{proof}
With the same notation as in the statement of the theorem, we have

$$
\int_{\Omega}\left(\sup_{t\in\mathbb{R}}\left|\langle{\delta_m, e^{-itJ_{\omega}}\delta_0}\rangle\right|\right)d\mu(\omega)
=a(m,0)\\
\leq\liminf_{L\to\infty}a_L(m,0)\\
\leq \liminf_{L\to\infty}\rho_L(m,0)\\$$
$$\leq \liminf_{L\to\infty}\frac{\sqrt{a_0a_{m-1}}}{a_{-L}a_{L-1}}\int_{\Sigma_0}\Big\langle{T_{E;1}^{(1)}\dots T_{E;m-1}^{(m-1)}S_{E;m}^{(m)}\dots S_{E;L-1}^{(L-1)}\phi_{L;E;L}^{(L-1)}, US_{E;0}^{(0)}\dots S_{E;-L+1}^{(-L+1)}\phi_{-L;E;-L}^{(-L)}}\Big\rangle_{L^2(\mathbb{R}, dx_1)}dE.\\$$
$$\leq \norm{a}_\infty\cdot \delta^{-2}\liminf_{L\to\infty}\int_{\Sigma_0}\norm{T_{E;1}^{(1)}\dots T_{E;m-1}^{(m-1)}S_{E;m}^{(m)}\dots S_{E;L-1}^{(L-1)}\phi_{L;E;L}^{(L-1)}}_2\norm{ US_{E;0}^{(0)}\dots S_{E;-L+1}^{(-L+1)}\phi_{-L;E;-L}^{(-L)}}_2dE.\\$$
$$=\norm{a}_\infty\cdot \delta^{-2}\liminf_{L\to \infty}\int_{\Sigma_0}\norm{T_{E;1}^{(1)}\dots T_{E;m-1}^{(m-1)}S_{E;m}^{(m)}\dots S_{E;L-1}^{(L-1)}\phi_{L;E;L}^{(L-1)}}_2\norm{ S_{E;0}^{(0)}\dots S_{E;-L+1}^{(-L+1)}\phi_{-L;E;-L}^{(-L)}}_2dE.\\$$
$$\leq \norm{a}_\infty\cdot \delta^{-2}\liminf_{L\to \infty}\int_{\Sigma_0}\norm{T_{E;1}^{(1)}\dots T_{E;m-1}^{(m-1)}}_{2,2}\norm{S_{E;m}^{(m)}}_{1,2}\norm{S_{E;m+1}^{(m+1)}}_{1,1}\dots \norm{S_{E;L-1}^{(L-1)}}_{1,1}\norm{\phi_{L;E;L}^{(L-1)}}_1\\$$
$$\times \norm{ S_{E;0}^{(0)}}_{1,2}\norm{S_{E;1}^{(1)}}_{1,1}\dots \norm{S_{E;-L+1}^{(-L+1)}}_{1,1}\norm{\phi_{-L;E;-L}^{(-L)}}_1dE.\\$$
$$\leq\norm{a}_\infty\cdot \delta^{-4}\liminf_{L\to\infty}\int_{\Sigma_0}\, q^{\frac{m-2}{2}}\sqrt {a_{m-1}\cdot a_{-1}}\norm{r}_\infty dE\\$$
$$\leq \norm{a}_\infty^2\cdot \delta^{-4}\norm{r}_\infty Leb(\Sigma_0)q^{\frac{m-2}{2}}\\$$
$$=C\cdot e^{-\gamma|m|},$$
where $\displaystyle C= \norm{a}_\infty^2\cdot \delta^{-4}\norm{r}_\infty Leb(\Sigma_0)q^{-1}$, and $\displaystyle\gamma=\frac{1}{2}\log \left(q^{-1}\right).$

\end{proof}

\subsection{Proof of Theorem \ref{thmm2.6}}

\begin{proof}
With the same notation as in the statement of the theorem, we have

$$
\int_{\Omega}\left(\sup_{t\in\mathbb{R}}\left|\langle{\delta_m, e^{-itJ_{\omega}}\delta_0}\rangle\right|\right)d\mu(\omega)
=a(m,0)\\
\leq\liminf_{L\to\infty}a_L(m,0)\\
\leq \liminf_{L\to\infty}\rho_L(m,0)\\$$
$$\leq \liminf_{L\to\infty}\frac{\sqrt{a_0a_{m-1}}}{a_{-L}a_{L-1}}\int_{\Sigma_0}\Big\langle{T_{E;1}^{(1)}\dots T_{E;m-1}^{(m-1)}S_{E;m}^{(m)}\dots S_{E;L-1}^{(L-1)}\phi_{L;E;L}^{(L-1)}, US_{E;0}^{(0)}\dots S_{E;-L+1}^{(-L+1)}\phi_{-L;E;-L}^{(-L)}}\Big\rangle_{L^2(\mathbb{R}, dx_1)}dE.\\$$
$$\leq \norm{a}_\infty\cdot \delta^{-2}\liminf_{L\to\infty}\int_{\Sigma_0}\norm{T_{E;1}^{(1)}\dots T_{E;m-1}^{(m-1)}S_{E;m}^{(m)}\dots S_{E;L-1}^{(L-1)}\phi_{L;E;L}^{(L-1)}}_2\norm{ US_{E;0}^{(0)}\dots S_{E;-L+1}^{(-L+1)}\phi_{-L;E;-L}^{(-L)}}_2dE.\\$$
$$=\norm{a}_\infty\cdot \delta^{-2}\liminf_{L\to \infty}\int_{\Sigma_0}\norm{T_{E;1}^{(1)}\dots T_{E;m-1}^{(m-1)}S_{E;m}^{(m)}\dots S_{E;L-1}^{(L-1)}\phi_{L;E;L}^{(L-1)}}_2\norm{ S_{E;0}^{(0)}\dots S_{E;-L+1}^{(-L+1)}\phi_{-L;E;-L}^{(-L)}}_2dE.\\$$
$$\leq \norm{a}_\infty\cdot \delta^{-2}\liminf_{L\to \infty}\int_{\Sigma_0}\norm{T_{E;1}^{(1)}\dots T_{E;m-1}^{(m-1)}}_{2,2}\norm{S_{E;m}^{(m)}}_{1,2}\norm{S_{E;m+1}^{(m+1)}}_{1,1}\dots \norm{S_{E;L-1}^{(L-1)}}_{1,1}\norm{\phi_{L;E;L}^{(L-1)}}_1\\$$
$$\times \norm{ S_{E;0}^{(0)}}_{1,2}\norm{S_{E;1}^{(1)}}_{1,1}\dots \norm{S_{E;-L+1}^{(-L+1)}}_{1,1}\norm{\phi_{-L;E;-L}^{(-L)}}_1dE.\\$$
$$\leq \norm{a}_\infty\delta^{-4}\liminf_{L\to\infty}\int_{\Sigma_0}A(1)\times\dots\times A\left(\floor*{\frac{m-1}{2}}\right)\sqrt{d_m^{-1}a_{m-1}}\sqrt{d_0^{-1}a_{-1}}\norm{r}_\infty dE$$
$$\leq \sqrt{d_0^{-1}}\norm{r}_\infty\norm{a}_\infty^2\delta^{-4}d_m^{-1/2}A(1)\times\dots\times A(k) $$
$$=\tilde C\, d_{m}^{-1/2}A(1)\times\dots\times A\left(\floor*{\frac{m-1}{2}}\right), $$
$$\leq \tilde C\, d_{m}^{-1/2}\exp\left(-\gamma' |k|^{1-2\zeta}\right).$$
$$\leq \tilde C \times C_1 |m|^{\zeta/2}\exp\left(-\gamma' \left|\floor*{\frac{m-1}{2}}\right|^{1-2\zeta}\right)$$
$$\leq C' |m|^{\zeta/2}\exp\left(-\gamma'' |m|^{1-2\zeta}\right)$$

where $\tilde C= \sqrt{d_0^{-1}}\norm{r}_\infty\norm{a}_\infty^2\delta^{-4}Leb(\Sigma_0),$ and $C'=\tilde C\times C_1.$

\end{proof}

\vskip 100pt
\noindent {\sc Acknowledgments:} The author would like to thank David Damanik for introducing him to this problem and also for his guidance throughout this project. The author is also grateful to Jake Fillman for his comments, which have improved the exposition of the paper.

\end{document}